\documentclass[11pt,a4paper]{article}

\pdfoutput=1

\usepackage[utf8]{inputenc}
\usepackage{amsmath,bm}
\usepackage{amsfonts}
\usepackage{amssymb}
\usepackage{upgreek}
\usepackage{fontenc}
\usepackage{comment}
\usepackage{graphicx}
\usepackage{accents}
\usepackage{trimclip}
\usepackage[title]{appendix}
\usepackage{caption}
\usepackage{subcaption} 
\usepackage{multirow}

\newcommand{\de}[0]{\mbox{ d}}
\newcommand {\der} [2] {\ds {\frac{\partial #1}{\partial #2}}}
\newcommand {\derl} [2] {(\ds {\partial #1/\partial #2})}

\newcommand {\ds}{\displaystyle}

\newcommand {\bvec} [1] {\pmb{\mathbf {#1}}}
\newcommand{\tsigma}{\pmb{\boldsymbol{\upsigma}}}
\newcommand {\tepsilon} {\pmb{\boldsymbol{\upvarepsilon}}}

\let \oldnabla \nabla
\renewcommand {\nabla} [0] {\bm{\oldnabla}}

\newcommand{\rn}[1]{%
  \textup{\uppercase\expandafter{\romannumeral#1}}%
}

\usepackage{stackengine}
\stackMath
\newcommand\tenq[2][1]{%
 \def\useanchorwidth{T}%
  \ifnum#1>1%
    \stackunder[0pt]{\tenq[\numexpr#1-1\relax]{#2}}{\scriptscriptstyle\sim}%
  \else%
    \stackunder[1pt]{#2}{\scriptscriptstyle\sim}%
  \fi%
}
\stackMath
\newcommand\mybar[2][1]{%
 \def\useanchorwidth{T}%
  \ifnum#1>1%
    \stackunder[0pt]{\mybar[\numexpr#1-1\relax]{#2}}{\scriptscriptstyle\textminus}%
  \else%
    \stackunder[1pt]{#2}{\scriptscriptstyle\textminus}%
  \fi%
}

\newcommand{\tensf}[1] 	{\underaccent \bar{{\bm{#1}}}}
\newcommand{\tenss}[1] 	{\underaccent \tilde{{\bm{#1}}}}

\newcommand{\tensff}[1] 	{\underaccent \tilde{\underaccent \tilde{{\bm{#1}}}}}

\begin{document}
\title{A micropolar isotropic plasticity formulation for non associated flow rule and softening featuring multiple classical yield criteria.   \newline Part II - FE integration and applications}
\author{Andrea Panteghini, Rocco Lagioia}
\maketitle

\section{Abstract}
A Finite Element procedure based on a  full implicit backward Euler predictor/corrector scheme for the Cosserat continuum is here presented. 
Since this   is based on invariants of the stress and couple stress tensors and on the spectral decomposition   of the former,  considerable benefits are achieved. 
The integration requires the solution of a single equation in a single unknown, which is a considerable improvement as compared to the system of seven or four equations required by  other approaches available in the literature for the Cauchy medium.
The scheme also allows for a very efficient treatment of the singularity which affects the apex of most of the existing yield and plastic potential surfaces.  Moreover,  no complications arise when some of the principal stresses coincide.
The algorithm has been implemented in a proprietary Finite Element program,  and used for the constitutive model proposed in \textit{part I} of this paper.
Numerical analyses have been conducted to simulate a biaxial compression test and a shallow strip footing resting on a Tresca,  Mohr-Coulomb,  Matsuoka-Nakai and Lade-Duncan soil.  
The benefits of the Cosserat  continuum over the Cauchy/Maxwell  medium are discussed considering   mesh refinement,  non-associated flow and softening behaviour.

\section{Introduction}
The  rotational degrees of freedom of the Cosserat continuum introduce a characteristic length which prevents the governing equations of static problems  from loosing  ellipticity  when localisations occur in materials  characterized by softening and/or non associated flow. 
The Cosserat continuum is then an effective tool for regularizing hill-posed  boundary value problems. Mesh sensitivity,  lack of convergence to the exact solution,  dependence of the thickness of the shear band on the mesh refinement can all  be effectively dealt with.  It can also be effective in preventing early crushing of numerical analyses  involving non-associated plastic flow. 

Constitutive models for the Cosserat continuum can also  be integrated using a full  implicit backward Euler scheme, although this is particularly time consuming.  However, since it is the only unconditionally stable method, even with very non linear yield and plastic potential surfaces such as those   adopted here,   its use is particularly recommended and it is indeed chosen  in this study. Moreover in combination with a consistent tangent operator, quadratic convergence in the structural Newton loop is achieved.

The classical approach to the implicit integration uses  the components of the strain tensor as the main variables, resulting in a system of seven equations in seven unknowns to be  iteratively solved.  
An alternative approach is that proposed by Tamagnini et al.  \cite{Tamagnini2002} and Borja et al. \cite{Borja2003} which is based on the observation that the tensorial derivative of an isotropic scalar valued function of that tensor is coaxial with the tensor itself.  The spectral decomposition of the stress and strain tensors is then exploited to use the principal strains as  the main variables,  thus reducing the  system of equations and unknown to four.  Complications  however occur when two of the principal stresses coincide (e.g.  Miehe \cite{Miehe1998}).

More recently a different approach has been presented by Panteghini and Lagioia \cite{PL2014} \cite{PL2018} which is based on the use of invariants as the main variables during the iterations.  This considerably improves the efficiency of the backward Euler scheme since  a single equation in a single unknown needs to be solved  to integrate the constitutive relationship. Moreover no difficulties arise when two or even three of the principal stresses coincide.
This approach presents also another considerable advantage.  The  De Souza Neto et al.  \cite{bibbia} method   can be used to deal with singularities in the yield and plastic potential surfaces, which results in an extremely stable integration algorithm.  
In fact, as reported by Panteghini and Lagioia \cite{PL2014} \cite{PL2018} numerical  analyses of a  shallow footing problem, with no lateral surcharge and no effective cohesion in the failure criterion can be carried out without any difficulty, provided that an associated flow is adopted.  The analysis of that  boundary value problem is usually considered not feasible. 
However, when non associated plasticity is considered even this approach is not sufficient to avoid very early crushing of the analyses. 

In this article the Panteghini and Lagioia \cite{PL2018}  scheme is extended to the Cosserat continuum. 
As shown in what follows the formulation  is  laborious, as it requires the evaluation of an elevated number of derivatives and of fourth order tensors.
However the rewards in terms of stability of the Finite Element program and speed of the analyses is extremely high.
The algorithm has been implemented in a  proprietary FE program and numerical analyses have been conducted to simulate a biaxial compression test on a specimen with a defect inclusion.

The effectiveness of the constitutive model presented in part I of this paper  and of the integration scheme  has been tested in extreme conditions. 
The very demanding footing problem previously described has been analysed using a high rate of softening, described by an exponential law,  and a high non-associativeness. 
Analyses were conducted using failure criteria which are significant in both research and engineering practice, such as the Mohr-Coulomb and the Matuoska-Nakai, showing the potentiality also in real world applications. 

\section{The backward Euler incremental initial value problem in terms of invariants}

\subsection{The general  case}
The incremental  elasto-plastic   initial value problem aims at the evaluation of all stresses, couple-stresses, strains and curvatures   at the end of a pseudo time interval $\left[ t_n, t_{n+1} \right]$ during which a given increment of  strains $\Delta \tenss{\gamma}$ and curvatures $\Delta \tenss{\chi}$ is applied and  all  quantities are known at time  $t_n$.  To keep the notation as simple as possible,  all quantities at time $t_{n+1}$ are indicated without the $n+1$ subscript whilst for  those at the beginning of the pseudo time interval the $n$ subscript is retained.

The integration algorithm presented in this section is an extension to the Cosserat continuum of that proposed by Panteghini and Lagioia \cite{PL2018} for the classical Cauchy medium.  It is a backward Euler predictor/corrector algorithm formulated in terms of invariants of stress, couple-stress, strain and curvature  tensors. At variance with other algorithms it requires the solution of a single equation in a single unknown, resulting in an extremely fast numerical integration.

The usual assumption is made that the  strains and curvatures can be decomposed into their recoverable, elastic,  and permanent, plastic,  components 
\begin{equation}
\begin{gathered}
\dot{\tenss{\gamma}}= \dot{\tenss{\gamma}}^e + \dot{\tenss{\gamma}}^p \\
\dot{\tenss{\chi}}= \dot{\tenss{\chi}}^e + \dot{\tenss{\chi}}^p
\end{gathered}
\end{equation}

The tentative assumption that  strain increments over the pseudo time interval are elastic leads to the classical definition of strain predictors
\begin{equation}
\begin{gathered}
\tenss{\gamma}^* =\tenss{\gamma}_n^e+\Delta \tenss{\gamma} \\
\tenss{\chi}^*= \tenss{\chi}_n^e+\Delta \tenss{\chi}
\end{gathered}
\end{equation}
which is here rewritten exploiting the decomposition of a general second order tensor into its symmetric and skew-symmetric components 
\begin{equation}
\begin{gathered}
\tenss{\varepsilon}^*=\tenss{\varepsilon}_n^e+ \Delta\tenss{\varepsilon} \\
\tenss{\omega}^*=\tenss{\omega}_n^e + \Delta\tenss{\omega} \\
\text{sym}\tenss{\chi}^*=\text{sym}\tenss{\chi}_n^e+\text{sym}\Delta\tenss{\chi} \\
\text{skw}\tenss{\chi}^*=\text{skw}\tenss{\chi}_n^e+\text{sym}\Delta\tenss{\chi}
\end{gathered}
\end{equation}
and also separating the spherical and the deviatoric parts (note that the skew-symmetric tensor is by definition deviatoric)
\begin{equation}
\begin{gathered}
\varepsilon_v^*=\text{tr} \tenss{\varepsilon}_n^e + \Delta {\varepsilon}_v \\
\tenss{e}^*=\tenss{e}_n^e+ \Delta\tenss{e} \\
\tenss{\omega}^*=\tenss{\omega}_n^e + \Delta\tenss{\omega} \\
\text{tr} \tenss{\chi}^*=\text{tr} \tenss{\chi}_n^e+  \text{tr}\Delta\tenss{\chi}\\
\text{sym}\tenss{g}^*=\text{sym}\tenss{g}_n^e+\text{sym}\Delta\tenss{g} \\
\text{skw}\tenss{g}^*=\text{skw}\tenss{g}_n^e+\text{skw}\Delta\tenss{g}
\end{gathered}
\end{equation}
where the standard definition of deviatoric tensors has been used
\begin{equation}
\begin{gathered}
\tenss{e}=\tenss{\varepsilon}-\frac{1}{3} \text{tr} \tenss{\varepsilon} \tenss{I}\\
\tenss{g}= \tenss{\chi} - \frac{1}{3} \text{tr} \tenss{\chi} \tenss{I}
\end{gathered}
\nonumber
\end{equation}
Since the elastic constitutive equations retrieved in \textit{Part I} of this paper are
\begin{equation}
\begin{gathered}
\tenss{\sigma}= K  \text{tr}{\tenss{\varepsilon}^e} \tenss{I} + 2 G \tenss{e}^e + 2 G_c \tenss{\omega}^e \\
\tenss{\mu}= K_c \text{tr}\tenss{\chi}^e \tenss{I} + 2 B \tenss{g}^e_{sym} + 2 B_c \tenss{g}^e_{skw}
\end{gathered}
\label{Eq_StressStrainSpherDev}
\end{equation}
the stress predictors can be evaluated
\begin{equation}
\begin{gathered}
p^*= K \varepsilon_v^* \\
\tenss{s}_{sym}^*= 2G \tenss{e}^* \\
\tenss{s}_{skw}^*= 2 G_c \tenss{\omega}^*\\
\text{tr}{\tenss{\mu}}^*=K_c \text{tr}{\tenss{\chi}}^*  \\
\tenss{m}_{sym}^*= 2 B \tenss{g}_{sym}^*\\
\tenss{m}_{skw}^*= 2 B_c \tenss{g}_{skw}^*
\end{gathered}
\end{equation}
The definition of the \textit{equivalent von Mises stress} formulated in the \textit{Part I} of this paper
\begin{equation}
\begin{split}
q= \left\lbrace \frac{3}{2} \left[ \tenss{s}_{sym} \colon \tenss{s}_{sym} + \frac{G}{G_c} \tenss{s}_{skw} \colon \tenss{s}_{skw} + \frac{G}{B} \tenss{m}_{sym} \colon \tenss{m}_{sym} +\frac{G}{B_c} \tenss{m}_{skw} \colon \tenss{m}_{skw} \right] \right.  \\
\left.   + \frac{2 G}{K_c} \frac{ \text{tr}^2\tenss{\mu}}{9}  \right\rbrace ^{\frac{1}{2}}
\end{split}
\label{Eq_vonMisesEqStressCosserat}
\end{equation}
can then be used to determine the 
 invariants of the stresses and couple stress predictor $\theta_s^*$ (Lode's angle),  $q_s^*$  and $q^*$ 
\begin{equation}
\begin{aligned}
q^*_s= \sqrt{\frac{3}{2} \tenss{s}_{sym}^* \colon \tenss{s}_{sym}^*},\;\;\;
\theta^*_s=\frac{1}{3} \arcsin \left[-\frac{27}{2} \frac{\text{det} \tenss{s}^*_{sym}}{q^*_s} \right]
\\
q^*= \left\lbrace \frac{3}{2} \left[ \tenss{s}_{sym}^* \colon \tenss{s}_{sym}^* + \frac{G}{G_c} \tenss{s}_{skw}^* \colon \tenss{s}_{skw}^* + \frac{G}{B} \tenss{m}_{sym}^* \colon \tenss{m}_{sym}^* +\frac{G}{B_c} \tenss{m}_{skw}^* \colon \tenss{m}_{skw}^*  \right.\right.  \\
\left.\left.  + \frac{2 G}{ K_c}  \frac{\text{tr}^2\tenss{\mu}^*}{9} \right]  \right\rbrace ^{\frac{1}{2}} \\
\end{aligned}
\end{equation}
where the subscript $s$ in the first two invariants refers to the symmetric part of the deviatoric stress tensor alone. 

All quantities for the evaluation of the  yield function 
\begin{equation}
f(p,q,\theta_s)= q \,  \Gamma \left(\theta_s \right) + M_c p- \sigma_0\left(\lambda \right)  
\label{Eq_YieldFunction}
\end{equation}
in the predictor conditions have been determined. 
The resulting value,  $f^*= f\left( p^*, q^*, \theta_s^*, \sigma_0(\lambda_n) \right)$  enables to assess whether the tentative assumption of a completely elastic strain/curvature increment was correct (i.e.  $f^* \le 0$),  hence resulting in all stress, couple stress,  strain and curvature  predictor quantities being the actual values  at the end of time $t_{n+1}$. 

If on the other hand  $f^*>0$ the strain/curvature increment applied during the  pseudo time interval $\left[ t_n, t_{n+1} \right]$ was elasto-plastic and a return to the appropriate yield surface must be performed to evaluate the correct elastic and plastic parts of that  increment
\begin{equation}
\begin{gathered}
\Delta \varepsilon_v=\Delta \varepsilon^e_v + \Delta \varepsilon_v^p \\
\Delta \tenss{e}=  \Delta \tenss{e}^e+ \Delta \tenss{e}^p\\
\Delta \tenss{\omega} = \Delta \tenss{\omega}^e + \Delta \tenss{\omega}^p \\
\text{tr} \Delta \tenss{\chi}= \text{tr} \Delta \tenss{\chi}^e +\text{tr} \Delta \tenss{\chi}^p \\
\text{sym}\Delta \tenss{g}= \text{sym}\Delta \tenss{g}^e+ \text{sym}\Delta \tenss{g}^p \\
\text{skw}\Delta \tenss{g}= \text{skw}\Delta \tenss{g}^e+ \text{skw}\Delta \tenss{g}^p
\end{gathered}
\label{Eq_StrainIncrElastPlastComp}
\end{equation}

The return algorithm requires the formulation  of all stress and couple stress at the end of the pseudo time interval  as a function of the respective predictor quantities, which is readily done with the exception of the symmetric component of the deviatoric stress tensor 
\begin{equation}
\begin{gathered}
p= p_n + K (\Delta \varepsilon_v- \Delta \varepsilon_v^p) = p_n^* - K \Delta \varepsilon_v^p\\
\tenss{s}_{skw}= \tenss{s}_{skw,n}+2 G_c (\Delta \tenss{\omega} - \Delta \tenss{\omega}^p)= \tenss{s}_{skw,n}^*- 2G_c \Delta \tenss{\omega}^p \\
\text{tr} \tenss{\mu}= \text{tr} \tenss{\mu}_n+  K_c (\text{tr} \Delta  \tenss{\chi} -  \text{tr} \Delta  \tenss{\chi}^p )= \text{tr}\tenss{\mu}_n^*  - K_c  \text{tr} \Delta \tenss{\chi}^p  \\
\tenss{m}_{sym}= \tenss{m}_{sym,n}+ 2B (\Delta \tenss{g}_{sym} - \Delta \tenss{g}_{sym}^p)= \tenss{m}_{sym,n}^* - 2B \Delta \tenss{g}_{sym}^p \\
\tenss{m}_{skw}= \tenss{m}_{skw,n}+2 B_c (\Delta \tenss{g}_{skw} - \Delta \tenss{g}_{skw}^p)= \tenss{m}_{skw,n}^* - 2B_c \Delta \tenss{g}_{skw}^p 
\end{gathered}
\label{Eq_nonstandards}
\end{equation}
The plastic components in the previous equation can be evaluated in the standard manner by deriving the plastic potential  
\begin{equation}
g(p,q,\theta_s)= q \,  \hat{\Gamma} \left(\theta_s \right) + M_c p 
\nonumber
\end{equation}
with respect to the associated stress quantities,  resulting in 
\begin{equation}
\begin{gathered}
\Delta \varepsilon_v^p=-\hat{M} \Delta\lambda \\
\Delta \tenss{\omega}^p= \frac{\partial g}{\partial \tenss{\sigma}_{sym}} \Delta \lambda = \frac{3}{2}\frac{G}{G_c} \frac{\hat{\Gamma}(\theta_s) \Delta \lambda} { q} \tenss{s}_{skw} \\
\text{tr} \Delta \tenss{\chi}^p=   \frac{\partial g}{\partial \frac{1}{3} \text{tr} \tenss{\mu}} \Delta \lambda = \frac{3}{2}\frac{2G}{K_c} \frac{\hat{\Gamma}(\theta_s) \Delta \lambda} { q}  \frac{\text{tr} \tenss{\mu}}{3} \\
\Delta \tenss{g}_{sym}^p= \frac{\partial g}{\partial \tenss{m}_{sym}} \Delta \lambda = \frac{3}{2}\frac{G}{B} \frac{\hat{\Gamma}(\theta_s) \Delta \lambda} { q} \tenss{m}_{sym} \\
\Delta \tenss{g}_{skw}^p= \frac{\partial g}{\partial \tenss{m}_{skw}} \Delta \lambda = \frac{3}{2}\frac{G}{B_c} \frac{\hat{\Gamma}(\theta_s) \Delta \lambda} { q} \tenss{m}_{skw} 
\end{gathered}
\end{equation}
where  $\Delta\lambda$ is the plastic multiplier. 
The stresses  and couple stress at time $t_{n+1}$ are hence  expressed  in terms of their respective predictor quantities
\begin{equation}
\begin{gathered}
p= p_{n}^*- K \hat{M} \Delta \lambda \\
\tenss{s}_{skw}=  \frac{ q}{q+3 G \, \hat{\Gamma}(\theta_s) \, \Delta \lambda} \: \tenss{s}_{skw}^*  \\
\tenss{\mu}= \frac{ q}{q+3 G \, \hat{\Gamma}(\theta_s) \, \Delta \lambda} \: \tenss{\mu}^*
\end{gathered}
\label{Eq_StresssFuncPredictors}
\end{equation}

This set of equations is not complete, since a similar expression would be  required for the symmetric component of the deviatoric stress tensor $ \tenss{s}_{sym}$ .
However since 
\begin{equation}
\Delta \tenss{e}^p=\der{g}{\tenss{s}_{sym}} \Delta \lambda=\left(\der{g}{p}\der{p}{\tenss{s}_{sym}}+\der{g}{q}\der{q}{\tenss{s}_{sym}}+\der{g}{\theta_s}\der{\theta_s}{\tenss{s}_{sym}}\right)\Delta \lambda
\label{Eq_DeltaTensEPlast}
\end{equation}
an explicit dependence on  $ \tenss{s}_{sym}$  on its predictor quantity cannot be obtained. 
 An alternative route can be followed, which was outlined by 
 Panteghini and Lagioia \cite{PL2018} 
who showed that a relationships can be formulated between the invariants of $\tenss{s}_{sym}$ and of $\tenss{s}_{sym}^*$ 
This is based on  geometrical considerations in the deviatoric plane and is  valid only if the elastic strain energy potential is not dependent on the third invariant, which results in the  deviatoric stress and strain tensors having the same Lode's angle.

\begin{figure}[tb]
\centering
\includegraphics[width=0.7\textwidth]{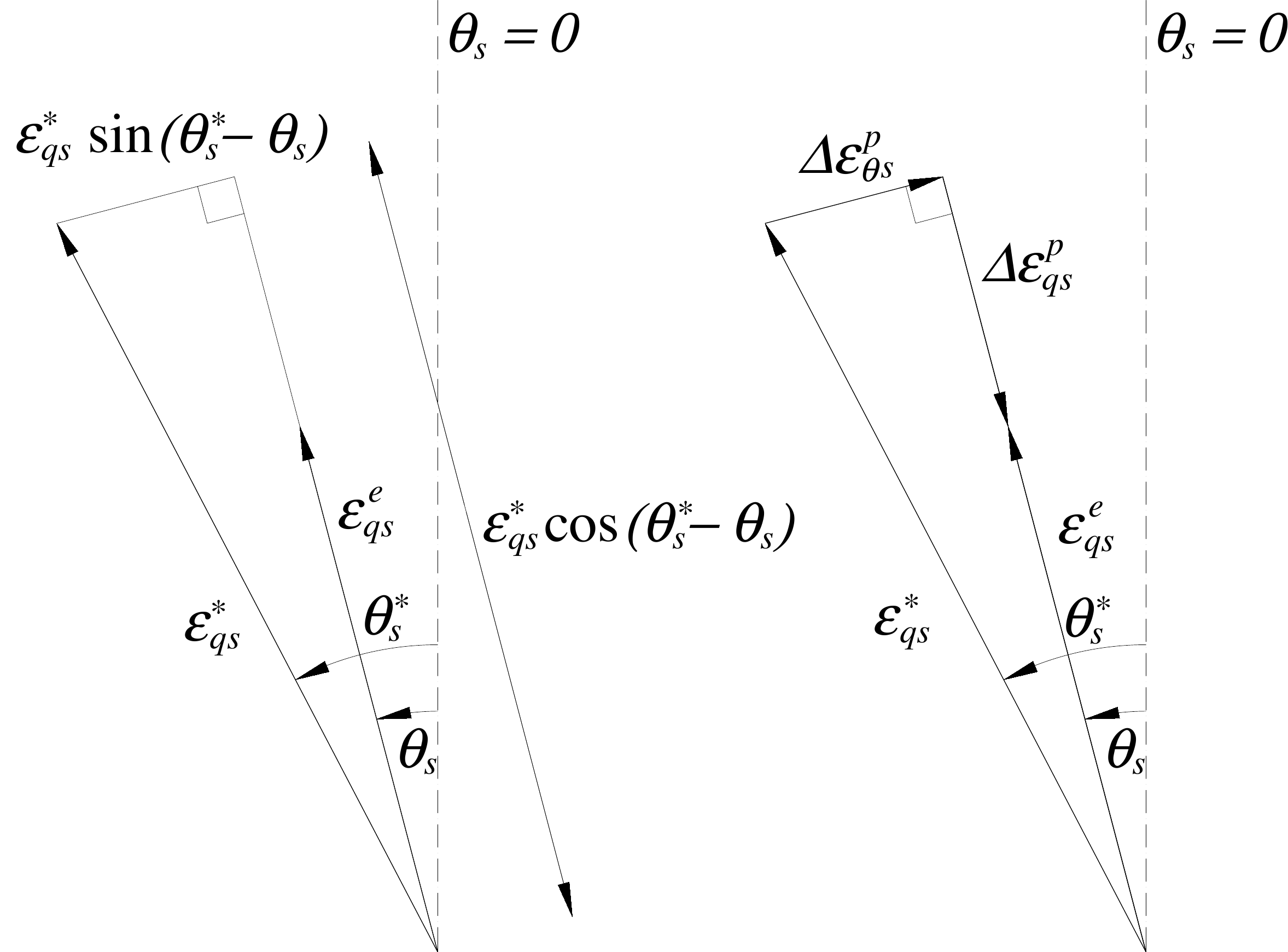}
\caption{Geometric representation in the deviatoric plane of the High-Westergaard principal symmetric strain space.  Elastic  predictor, plastic and final  elastic strain (modified from Panteghini and Lagioia   \cite{PL2018})}.
\label{fig:legametrig}
\end{figure}

Fig.\ref{fig:legametrig}  shows a portion of the deviatoric plane in the Haigh-Westergaard principal symmetric strain space. 
In that figure $\varepsilon_{qs}$ indicates the second invariant of the symmetric deviatoric strain tensor which is energetically conjugated to the equivalent von Mises stress $q_s$, and  $^*$  and $^e$ superscripts are as usual used  to distinguish the predictor and the elastic deviatoric strains  associated to the pseudo time interval $\left[t_n, t_{n+1} \right]$
\begin{equation}
\varepsilon_{qs}^e = \ds \sqrt{ \frac{2}{3} \tenss{e}^e: \tenss{e}^e}, \;\;\; \varepsilon_{qs}^* = \ds \sqrt{ \frac{2}{3} \tenss{e}^*: \tenss{e}^*}, \;\;\; q_s=\ds \sqrt{ \frac{3}{2} \tenss{s}_{sym} : \tenss{s}_{sym}}
\nonumber
\end{equation}

The figure shows that the actual elastic strain at the end of an  elasto-plastic strain increment is obviously different from the predictor quantity, and that  they are associated to different Lode's angle. The right hand-side diagram shows that the difference between those strains represents the plastic strain increment during that time interval and the radial and circumferential  components are highlighted.

The figure shows that the radial plastic component is
\begin{equation}
\Delta \varepsilon_{qs}^p = \varepsilon_{qs}^{*} \cos \left(\theta_s^*-\theta_s \right) -\varepsilon_{qs}                                                                                                                                                                                                                                                                                                                                                                                                                                                                                                                                                                                                                                                                                                                                                                                                                                                                                                                                                                                                                                                                                                                                                                                                                                                                                                                                                                                                                                                                            ^{e}
\label{eq:trigeqs1}
\end{equation}
and by applying the elastic constitutive law 
\begin{equation}
q_s = 3 G  \varepsilon_{qs}^e
\label{eq:legameqs}
\end{equation}
to Eq. \eqref{eq:trigeqs1}  and solving for $q_s$ one obtains
\begin{equation}
q_s=q_s^*\cos \left( \theta_s^*-\theta_s\right)-3 G \Delta \varepsilon_{qs}^p
\label{eq:tmpqs}
\end{equation}
Finally substituting   the plastic strain increment definition
\begin{equation}
\Delta \varepsilon_{qs}^p=\der{g}{q_s}\Delta \lambda=\frac{q_s \hat \Gamma\left( \theta_s\right) \Delta \lambda}{q}
\end{equation}
into Eq. \eqref{eq:tmpqs}, and collecting  $q_s$ yields
\begin{equation}
q_s=\left(\frac{ q}{ q + 3 G  \hat \Gamma\left(\theta_s\right) \Delta \lambda} \right) q_s^* \cos \left(\theta_s^*-\theta_s \right)
\label{Eq_qsFuncPredictors}
\end{equation}
which integrates  Eqs.\eqref{Eq_StresssFuncPredictors} in the definition of stress quantities at the end of the pseudo time interval as a function of the respective predictor.
If this equation is further introduced into Eq. \eqref{Eq_vonMisesEqStressCosserat}  
one obtains
\begin{equation}
\begin{gathered}
q=\left\{q_s^2\cos^2 \left(\theta_s^*-\theta_s \right)+\frac{3}{2}\left[ \left(\frac{G }{G_c} \right)\tenss{ s}^*_{skw}:\tenss{ s}^*_{skw}+\left(\frac{G}{B} \right)\tenss{ m}^*_{sym}:\tenss{ m}^*_{sym}\right.\right.\\
\left.\left.+\left(\frac{G }{B_c} \right)\tenss{ m}^*_{skw}:\tenss{ m}^*_{skw}  + \frac{2 G}{ K_c} \frac{ \text{tr}^2\tenss{\mu}^*}{9} \right]\right\}^{\frac{1}{2}}- 3 G \hat \Gamma\left(\theta_s\right) \Delta \lambda
\end{gathered}
\label{Eq_VonMisesProva}
\end{equation}
which can be furthermore simplified by observing that
\begin{equation}
\begin{gathered}
\frac{3}{2}\left[  \left(\frac{G }{G_c} \right)\tenss{ s}^*_{skw}:\tenss{ s}^*_{skw}+\left(\frac{G}{B} \right)\tenss{ m}^*_{sym}:\tenss{ m}^*_{sym} +\left(\frac{G }{B_c} \right)\tenss{ m}^*_{skw}:\tenss{ m}^*_{skw}  + \right. \\
\left.  \frac{2 G}{ K_c} \frac{ \text{tr}^2\tenss{\mu}^*}{9} \right]= 
{q^*}^2-{q_s^*}^2
\end{gathered}
\nonumber
\end{equation}
so that  Eq.  \ref{Eq_VonMisesProva} becomes
\begin{equation}
q=r-3 G  \hat \Gamma\left(\theta_s\right) \Delta \lambda
\label{Eq_qFuncRLambdaTheta}
\end{equation}
where
\begin{equation}
r=\sqrt{{q^*}^2-{q_s^*}^2\sin^2\left(\theta^*_s-\theta_s \right)}
\label{Eq_r}
\end{equation}
In conclusion, the developments  so far presented,  result in 
\begin{equation}
\begin{gathered}
p= p_{n}^*+ K \hat{M} \Delta \lambda \\
q_s=\left(\frac{q}{r}\right)  q_s^*\cos \left(\theta^*_s-\theta_s \right) \\
\tenss{ s}_{skw} = \left(\frac{q}{r}\right) \tenss{s}_{skw}^*\\
\tenss{\mu}=\left(\frac{q}{r}\right)  \tenss{\mu}^*\\
\end{gathered}
\label{Eq_StressesFuncPredicLambdaTheta}
\end{equation}
where the definition of $r$ has been substituted in Eqs. \eqref{Eq_StresssFuncPredictors}. 
 The meaning of the last set of equations is that the stress invariants $p$ and $q$  at the end of the pseudo time interval $\left[ t_n, t_{n+1}] \right]$ are a function of the known predictor invariants $p^*$, $q^*$, $q_s^*$ and $\theta_s^*$ and of two unknowns,  namely  the Lode's angle $\theta_s$ and the plastic multiplier $\Delta \lambda$ at the end of time interval  $t_{n+1}$.
However, as shown by Panteghini and Lagioia \cite{PL2014} \cite{PL2018},   these two unknown are not independent as $\Delta \lambda$ is a function of  $\theta_s$ 
\begin{equation}
\Delta \lambda= \frac{{q_s^*}^2 \frac{1}{2}\sin \left(2\theta_s^*-2\theta_s \right) }{3 G \, r \,  \hat \Gamma' \left( \theta_s\right) } 
\label{Eq_DeltaLamdaTheta}
\end{equation}
as can be easily shown  on the basis of geometrical reasoning from Fig. \ref{fig:legametrig}
\begin{equation}
\Delta \varepsilon_{\theta s}^p=\varepsilon_{qs}^{*} \sin \left(\theta_s^*-\theta_s \right) 
\label{Eq_EpsPlasticThetaGeometry}
\end{equation}
and substituting $\varepsilon_{qs}^*$ with $q^*/(3 G)$ using the linear elastic relationship and $\Delta \varepsilon_{qs}^p$  with its definition using the plastic potential
\begin{equation}
\Delta \varepsilon_{\theta s}^p=\frac{1}{q_s}\der{g}{\theta_s}\Delta \lambda=\frac{r \,  \hat \Gamma' \left( \theta_s\right) \Delta \lambda}{q_s^*\cos \left(\theta^*_s-\theta_s \right)}
\label{Eq_EpsPlasticTheta}
\end{equation}

The whole  framework on which the integration algorithm is based  has been presented.  The important results is that  all stress quantities at the end of the pseudo time interval $\left[ t_n, t_{n+1}\right]$ are a function of the predictor quantities, which are known, and of the plastic multiplier $\Delta \lambda$ alone. Hence the numerical integration of the constitutive law is obtained by iterating on a single equation in one unknown. This brings considerable advantages in terms of machine run-time as compared to classical approaches where a system of seven equations needs to be solved and also to approaches based on principal stresses  where the number of equations in the system is reduced to four.

\subsection{A particular case}
The general mathematical formulation of the return  algorithm presented in the previous subsection shows that  the converged Lode's angle  $\theta_s$  at time $t_{n+1}$  of the symmetric part of the stress tensor    differs from its predictor $\theta_s^*$.   
Eq. \eqref{Eq_DeltaLamdaTheta} is then exploited  to reduce the number of unknowns from two to one. 

\begin{figure}[tb]
\centering
\includegraphics[height=0.5\textwidth]{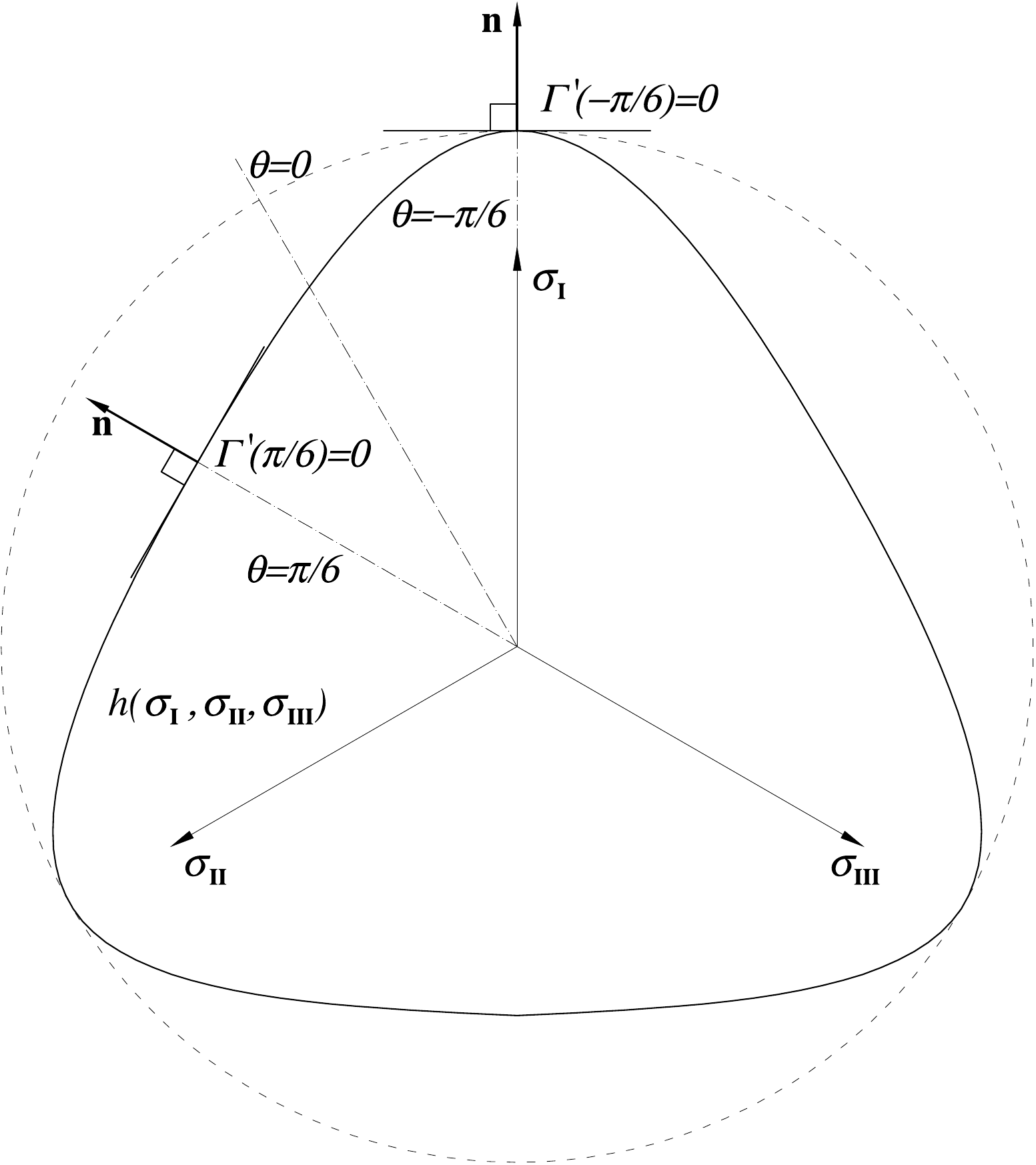}
\caption{Deviatoric section of the plastic potential surface (from Panteghini and Lagioia \cite{PL2018})}
\label{Fig_PlasticNormal}
\end{figure}

However a  particular case also occurs, which allows a significant  simplification of the integration algorithm.
Fig.\ref{Fig_PlasticNormal} shows the plastic potential section in the  deviatoric plane of  the Haigh-Westergaard symmetrical principal stress space. 
For symmetry reasons the derivative of the plastic potential with respect to the Lode's angle  $ \hat \Gamma'(\theta_s)$  vanishes when  $\theta_s=-\frac{\pi}{6}$ and $\theta_s=\frac{\pi}{6}$.
Since a particular  instance of an implicit integration scheme is here presented, hence  plastic strains are evaluated using the derivatives of the plastic potential   at the end of pseudo-time interval,   the converged  circumferential plastic strain $\Delta \varepsilon_{\theta s}^p$ also   vanishes, as shown  by Eq. \eqref{Eq_EpsPlasticTheta}.

Specific shapes  of the plastic potential in the deviatoric plane,  e.g.  those provided by  classical yield/failure criteria,   might result in additional Lode's angle values where that derivative vanishes.  As an example this does not occur in the case of the  Matsuoka-Nakai criterion,  whilst an additional $\frac{\pi}{6}\geq\theta_s\geq-\frac{\pi}{6}$ exists in the case of the Mohr-Coulomb criterion. 
 \footnote{The original Mohr-Coulomb criterion presents cusps in the deviatoric plane,  however smooth versions of that criterion are part of the General Classical yield function retrieved by Lagioia and Panteghini \cite{LP2016} also with the defining parameters derived for that function by Lester and Sloan \cite{Lester2017}.}
Another interesting instance is that provided by the von Mises and Drucker-Prager criteria. These  are characterized by a circular deviatoric section, hence   $ \hat \Gamma'(\theta_s)$ is identically nil.

When   $\Delta \varepsilon_{\theta s}^p$ vanishes  no plastic  correction of the strain predictor  $\varepsilon_q^*$ is possible in the circumferential direction of the deviatoric plane, hence $\theta_s^*$ must coincide with the converged $\theta_s$ (Fig. \ref{fig:legametrig}). This  also emerges from Eq. \eqref{Eq_EpsPlasticThetaGeometry}.
 A simplified radial return algorithm can then be used.

It should be noted that 
 not only is $\theta_s^*=\theta_s$  when $ \hat \Gamma'(\theta_s)=0$, but also the reverse is true,  i.e. when $ \hat \Gamma'(\theta_s^*)=0$ then $\theta_s=\theta_s^*$. 
 In fact if  $\theta_s^*$  is such that  $ \hat \Gamma'(\theta_s^*)=0$  then the  converged solution cannot be other than  $\theta_s=\theta_s^*$,  otherwise two stress predictor characterized with the same $\theta_s^*$ and   $ \hat \Gamma'(\theta_s^*)=0$  could be associated to different  values of $\theta_s$, one of which with  $ \hat \Gamma'(\theta_s)\neq0$. This would not be  not compatible with the fact that  the set of algebraic solving equations is of implicit type and uses  stress predictors as principal variables.

The main implication of this argument is  that   given $\theta_s^*$ it is \textit{a priori} known whether a generic or a radial return algorithm can be used,  depending only on whether or not $ \hat \Gamma'(\theta_s^*)=0$ .
Moreover, when the radial return algorithm applies,   since the converged $\theta_s$ is  known already at the beginning of the pseud-time interval,  it is not anymore an unknow in  Eqs \eqref{Eq_StressesFuncPredicLambdaTheta}, which only depends  on $\Delta \lambda$ (i.e. there is no need of  Eq.\eqref{Eq_DeltaLamdaTheta}).


\subsection{Return algorithms}
The arguments presented in the previous subsection indicate that a general and a radial return algorithms are required. 
However  many constitutive models, e.g. those used in perfect plasticity, adopt classical yield/failure criteria to define the yield and the plastic potential surfaces ( non associated plasticity typically  being achieved by choosing different sets of defining parameters).
Such criteria are characterized by  linear meridional sections which generate  a cusp at the apex of the surface toward the origin of the stress space.  
In  the case of the plastic potential this results  in an undefined  plastic flow normal,  which is particularly problematic when  backward Euler implicit integration schemes are adopted.  
However this issue  can be very effectively handled using  the approach proposed by de Souza Neto et al. \cite{bibbia}, which then requires the formulation of a third return algorithm. 

The cusp  at the apex is associated to the lack of continuity of the first derivative of its defining function and has been a major topic of discussion and research, particularly in the Soil Mechanics community.  It has very often been considered an unwelcome feature of numerous  functions and many attempts have been made to smoothen their apex  (e.g.   Abbo and Sloan \cite{Abbo1995},  Panteghini and Lagioia \cite{PL2014a}).
\begin{figure}[tb]
\centering
\includegraphics[width=\textwidth]{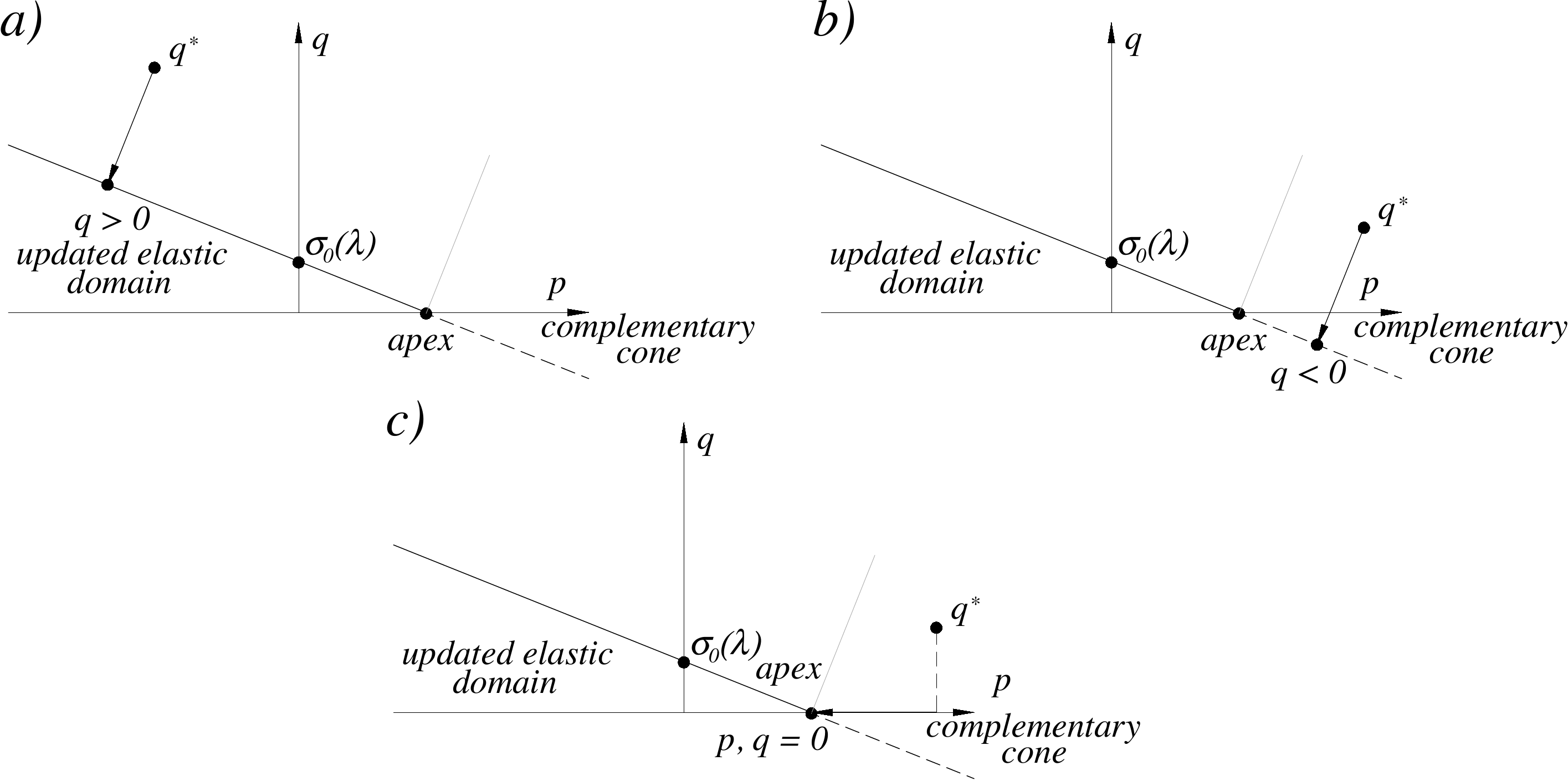}
\caption{Selection of the appropriate return mapping algorithm. Modified from \cite{bibbia}}
\label{fig:DP}
\end{figure}

However Panteghini and Lagioia \cite{PL2014}  performed  Finite Element analyses of a shallow foundation problem using a classical  criteria to define the yield and plastic potential surfaces  characterized by the apex singularity and compared the results of  either modifying the original formulation to provide smoothness or adopting  the   de Souza Neto et al.  \cite{bibbia} approach. The analyses showed that the latter approach is  considerably more efficient, resulting in reduced machine run-time.

Since the constitutive model adopted in this study also  uses classical yield/failure criteria together with the de Souza Neto et al. approach, three different return algorithms are required 
The first two  apply the mathematical developments described in the previous subsections and are mutually exclusive,  the one to be used being  know  \textit{a priori} on the basis of the predictor  Lode's angle.
 The first  algorithm is valid when $\hat \Gamma' \left( \theta^*_s\right)= 0$ which corresponds  for plastic potentials with  non-circular deviatoric sections  to  predictor Lode's angles $\theta_s^*=\pm\frac{\pi}{6}$, whilst for other   (e.g. the Mohr-Coulomb failure criterion)  also to  an additional intermediate value which only depends on the angle of shearing resistance $\phi$.   The second algorithm is valid for the general case of $\hat \Gamma' \left( \theta^*_s\right) \neq 0$,  i.e.   $\theta_s^*\neq\pm\frac{\pi}{6}$ (and for the intermediate Lode's angle for the Mohr-Coulomb criterion).

To deal with the cusp at the apex of the surface a control of the validity of the converged solution is needed, which in case of violation requires the use of the third algorithm. The validity check is based on the simple observation that the \textit{equivalent von Mises} stress $q$ is by definition a positive quantity.  As shown in Fig. \ref{fig:DP} for a generic meridional section of the  surface, given a generic predictor $q^*$ the return algorithm either results in a non-negative or in a negative value of $q$. In the former case the solution is clearly correct and can be accepted, whereas the case of a negative $q$ is  theoretically incorrect.
Since the solution is unique and the positive $q$ is excluded because it is characterized by a uniquely  defined plastic normal, then the only remaining valid possibility is a return to an isotropic stress state which is associated to an undefined plastic flow normal.   Hence a  third return algorithm needs to be formulated to bring the stress state back to the apex of the yield surface.


\paragraph{First return algorithm: $\hat \Gamma' \left( \theta^*_s\right)= 0$. }  \mbox{} \\
Considerable simplifications occurs, as  Eqs.  \eqref{Eq_EpsPlasticTheta}, \eqref{Eq_EpsPlasticThetaGeometry}  and \eqref{Eq_r}  indicate that  $\Delta \varepsilon_{\theta s}^p=0$,    $\theta_s=\theta^*_s$ and  $r\equiv q^*$. 
The stress invariants at the end of the increment can then be computed by solving a single non-linear scalar equation in the unknown $\Delta \lambda$ according to the step sequence:
\begin{itemize}
\item choose a starting trial value for $\Delta \lambda$; 
\item evaluate  $\sigma_0(\Delta \lambda)$ using the hardening/softening law;
\item compute $p$ using the first of Eqs.  \eqref{Eq_StressesFuncPredicLambdaTheta};
\item compute q from Eq. \eqref{Eq_qFuncRLambdaTheta},
\item   iterate on  $\Delta \lambda$ until  $f= q \Gamma \left(\theta^*_s \right) +M p -\sigma_0\left(\lambda \right)  = 0$ 
\footnote{
It should be noted that if a linear isotropic hardening/softening law is adopted 
\begin{equation}
\sigma_0(\lambda)=\bar \sigma_0 + h \left(\lambda_n+\Delta \lambda\right)
\label{eq:c0lin}
\end{equation}
then a close form solution for  $\Delta \lambda$ can be obtained 
\begin{equation}
\Delta \lambda=\frac{q^* \Gamma \left(\theta^*_s \right) + M p^* -\bar \sigma_0 -h \lambda_n }{   \Gamma \left(\theta^*_s \right) \hat \Gamma \left(\theta^*_s \right) 3G +  \left(M \hat M K  +h \right)}
\label{Eq_DeltaLambdaClosedForm}
\end{equation}
};
 \item check the validity of the converged solution: If $q\geq0$ the solution is valid, otherwise the third algorithm needs to be run to return to the apex of the yield surface.
\end{itemize} 

If the converged solution is valid,  all stress tensors at time $t_{n+1}$ can be easily evaluated,  in terms of the stress invariants $p$, $q$ and $\theta_s$ and of the predictor quantities using Eqs. \eqref{Eq_StressesFuncPredicLambdaTheta} and remembering that since   $\theta_s=\theta^*_s$  then   $r\equiv q^*$ and the stress and couple stress tensors can be updated using the following equations
 \footnote{The symmetric component of the deviatoric stress tensor at time $t_{n_1}$ is also proportional to $q/q^*$. 
Since in the first return algorithm  $\hat{\Gamma}'(\theta_s)=0$,   only  the second term on the right hand side of  Eq. \eqref{Eq_DeltaTensEPlast} in other than zero 
\begin{equation}
\Delta \tenss{e}^p=\der{g}{q}\der{q}{\tenss{s}_{sym}} \Delta \lambda=\frac{3}{2}\frac{\hat \Gamma \left(\theta^*_s \right) \Delta \lambda}{q} \tenss{s}_{sym}
\nonumber
\end{equation}
Morover since
\begin{equation}
\Delta \tenss{e}=\Delta \tenss{e}^e+\Delta \tenss{e}^p
\nonumber
\end{equation}
and
\begin{equation}
\tenss{s}_{sym} =\tenss{s}_{sym,n}+2 G\left(\Delta \tenss{e}-\Delta \tenss{e}^p\right)
\nonumber
\end{equation}
one obtains
\begin{equation}
\tenss{s}_{sym}=\frac{q}{q + 3 G \hat \Gamma\left(\theta^*_s\right) \Delta \lambda}  \tenss{s}_{sym}^*=\frac{q}{q^*}  \tenss{s}_{sym}^*
\label{eq:ssdp}
\end{equation}
where 
\begin{equation}
\tenss{s}_{sym}^* = \tenss{s}_{sym, n}+2 G  \Delta \tenss{e}
\nonumber
\end{equation}
}

\begin{equation}
\begin{gathered}
p=p^*_n+ K \hat{M} \,  \Delta \lambda\\
\tenss{s}_{sym}= \frac{q}{q^*} \tenss{s}_{sym}^*\\
\tenss{s}_{skw}= \frac{q}{q^*} \tenss{s}_{skw}^*\\
\tenss{\mu}=\left(\frac{q}{q^*}\right)  \tenss{\mu}^*\\
\end{gathered}
\label{Eq_StressTensorsFirstAlghorithm}
\end{equation}
The symmetric part of the stress tensor is then 
\begin{equation}
\tenss{\sigma}_{sym}= p \tenss{I} + \tenss{s}_{sym}=  p \tenss{I} + \frac{q}{q^* }\tenss{s}_{sym}^* 
\label{Eq_StressTensorsCosseratSymFirstAlgo}
\end{equation}
whilst the Cosserat stress is
\begin{equation}
\tenss{\sigma}= p \tenss{I} + \tenss{s}_{sym}+\tenss{s}_{skw} \\
\label{Eq_StressTensorsCosseratComponents}
\end{equation}

\paragraph{Second return algorithm: $\hat \Gamma' \left( \theta^*_s\right)\ne 0$.} \mbox{} \\
This is the general return algorithm for plastic potentials with a non-circular deviatoric section.  
In this case the  Lode's angle  $\theta_s$  at the end of the pseudo time interval is different from its predictor  $\theta_s^*$. 
As in the previous situation,  stresses and couple stresses at time $t_{n+1}$ are  obtained by solving a single non-linear scalar equation, but iteration are now performed on  $\theta_s$ (for example by employing a dumped Newton's method,  restricting the solution in the range  $-\pi/6\le \theta_s \le \pi/6$) rather than on the plastic multiplier $\Delta \lambda$
according to the step sequence: 
\begin{itemize}
\item choose a starting trial value for $\theta_s$;
\item evaluate  $\Delta \lambda$  using Eq. \eqref{Eq_DeltaLamdaTheta} as a function  of $\theta_s$ (and of the predictors  of the stress invariants);
\item evaluate  $\sigma_0(\Delta \lambda)$ using the hardening/softening law;
\item evaluate $p$ using the first of Eqs. \eqref{Eq_StressesFuncPredicLambdaTheta};
\item evaluate  $r$ and $q$ using Eqs. \eqref{Eq_r} and  \eqref{Eq_qFuncRLambdaTheta} ;
\item  iterate on  $\theta_s$ until  $f= q \Gamma \left(\theta^*_s \right) +M p -\sigma_0\left(\lambda \right)  = 0$;
\ \item check the validity of the converged solution: If $q\geq0$ the solution is valid, otherwise the third algorithm needs to be run to return to the apex of the yield surface.
\end{itemize}

If the converged solution is valid, the stresses at time $t_{n+1}$ need to be evaluated following the following procedure. 
\begin{itemize}
\item evaluate $q_s$ from the second of  Eqs.\eqref{Eq_StressesFuncPredicLambdaTheta};
\item evaluate the ordered principal stresses of the symmetric part of the stress tensor $\tenss{\sigma}_{sym}$ at time $t_{n+1}$
\begin{equation}
\begin{gathered}
\sigma_\rn{1}=p+\frac{2}{3} q_s \sin \left( \theta_s+\frac{2}{3} \pi \right)\\
\sigma_\rn{2}=p+\frac{2}{3} q_s \sin \left( \theta_s\right)\\
\sigma_\rn{3}=p+\frac{2}{3} q_s \sin \left( \theta_s-\frac{2}{3} \pi \right)\\
\end{gathered}
\end{equation}
where $\sigma_{I}\geq\sigma_{II}\geq\sigma_{III}$ \footnote{Note that since $\theta_s\neq \pm \frac{\pi}{6}$  these are strict inequalities.};
\item evaluate  $\tenss{\sigma}_{sym}$ by applying the spectral theorem
\begin{equation}
\tenss{\sigma}_{sym}= \sum_{i=\rn{1}, \rn{2}, \rn{3}} \sigma_i \tenss{b}_i^*
\label{Eq_spectral}
\end{equation}
where $\tenss{b}_i^*= \tensf{n}_i^* \otimes \tenss{n}_i^*$ and $\tensf{n}_i^*$ for the hypothesis of isotropy are the principal directions of  $\tenss{\varepsilon}^*$, which are know quantities for each pseudo-time step. 
\footnote{
It should be noted that the evaluation of the eigenvectors of the strain predictor $\tenss{\varepsilon}^*$ is actually not necessary,  as,  following Miehe (1998), \cite{Miehe1998} $\tenss{b}_i^*$ can be directly evaluated 
\begin{equation}
\tenss{b}_i^* \equiv \der{\varepsilon_i^*}{\tenss{\varepsilon}}=\der{\varepsilon_i^*}{\tenss{\varepsilon}^*}
\nonumber
\end{equation} 
where $\varepsilon_i^*$ is the $i-$th ordered principal component of $\tenss{\varepsilon}^*$
}

\item evaluate the remaining stress tensors using Eqs. \eqref{Eq_StressesFuncPredicLambdaTheta}; 
\item evaluate the Cosserat stress  tensor using Eq. \eqref{Eq_StressTensorsCosseratComponents}; 
\end{itemize}

\paragraph{Third return algorithm: $q<0$} \mbox{} \\
If either of the previous return algorithms converge to   $q<0$ the solution  is not theoretically acceptable and a return to the apex is performed. 
The yield condition then reduces to
\begin{equation}
M p - \sigma_0 (\lambda) = 0
\end{equation}
and after accountig for the  first of Eqs. \eqref{Eq_StressesFuncPredicLambdaTheta} a single non-linear scalar equation in $\Delta \lambda$ is obtained
\footnote{
Let note that $\Delta \lambda$ can be computed analytically from Eq. \eqref{eq:fapex}  if the linear strain hardening-softening law of Eq. \eqref{Eq_DeltaLambdaClosedForm} is employed.
In this case it results:
\begin{equation}
\Delta \lambda = \frac{M p^*  - \bar \sigma_0- h \lambda_n}{K M \hat M + h}
\end{equation}
}
\begin{equation}
M p^* -K M \hat M \Delta \lambda  - \sigma_0 (\lambda) = 0
\label{eq:fapex}
\end{equation}
whose iterative solution provides the value of $\Delta\lambda$ at time $t_{n+1}$, whilst all stress tensor are obtained from Eqs. \eqref{Eq_StressesFuncPredicLambdaTheta} 
\begin{equation}
\begin{gathered} \\
\tenss{\sigma}= \left(p^*- K \hat M \Delta \lambda \right) \tenss{I}\\
\tenss{s}_{sym}=\tenss{s}_{skw} = \tenss{\mu}=  \tenss{0}\\
\end{gathered}
\end{equation}

\section{Consistent jacobian matrices}
The consistent tangent operator is required to achieve quadratic convergence in the structural Newton loop (e.g Nagtegaal \cite{Nagtegaal1982},  Simo and Taylor, \cite{Simo1985}). 
 In algorithms for classical Cauchy continuum it is a single  fourth order tensor obtained by deriving   with respect to the strains the constitutive function $\hat {\tenss{\sigma}} (\tenss{\varepsilon})$ which provides the stresses as a function of the  strains   at time $t_{n+1}$
\begin{equation}
\der{\hat{\tenss{\sigma}}(\tenss{\varepsilon})}{\tenss{\varepsilon}}
\end{equation}
For the  Cosserat continuum  the number of the consistent operators increases together with the number of  stress tensors.

Since  the approach presented in this paper to integrate the constitutive law exploits  the spectral decomposition of the symmetric part of the  deviatoric  stress tensor,  the computation of the consistent operator involves the  derivation of the  eigenvalues and eigenvectors (or alternatively of the   bases) of that tensor.
This is not a straight forward task for the special occurrences  of  two  (i.e.  $\theta_s=\pm \pi/6$) or three (i.e. $q_s=0$)  coinciding  principal stresses, when some of the derivatives become singular  (e.g.  Miehe \cite{Miehe1998},  de Souza Neto et al. \cite{bibbia}). 

However, such difficulties can be easily overcome by formulating different consistent tangent operators for each of the return algorithms described in the previous section.
When the second return algorithm  applies (hence the solution is valid and $q_s>0$) all eigenvalues are distinct since $\theta_s \neq \pm \frac{\pi}{6}$.
The derivatives of the spectral decomposition can be evaluated without problems following the approach presented by Miehe \cite{Miehe1998}. 

On the other hand when the first or the third return  algorithms applies  two or three of  principal components, respectively,  of $\tenss{\sigma}_{sym}$ coincide and singularities would arise in the derivation of the basis  $\tenss{b}_i^*=\tenss{n}_i ^* \otimes \tenss{n}_i^*$.
However,  as discussed in the previous section , no  spectral decomposition is required in both cases, since when the first return algorithm applies, stresses are proportional to their predictors, whilst when the third algorithm applies the symmetric tensor is spherical. 

All necessary derivatives are provided in the appendix.

\section{Biaxial compression tests}
\begin{figure}[tb]
\centering
\includegraphics[height=0.5\textwidth]{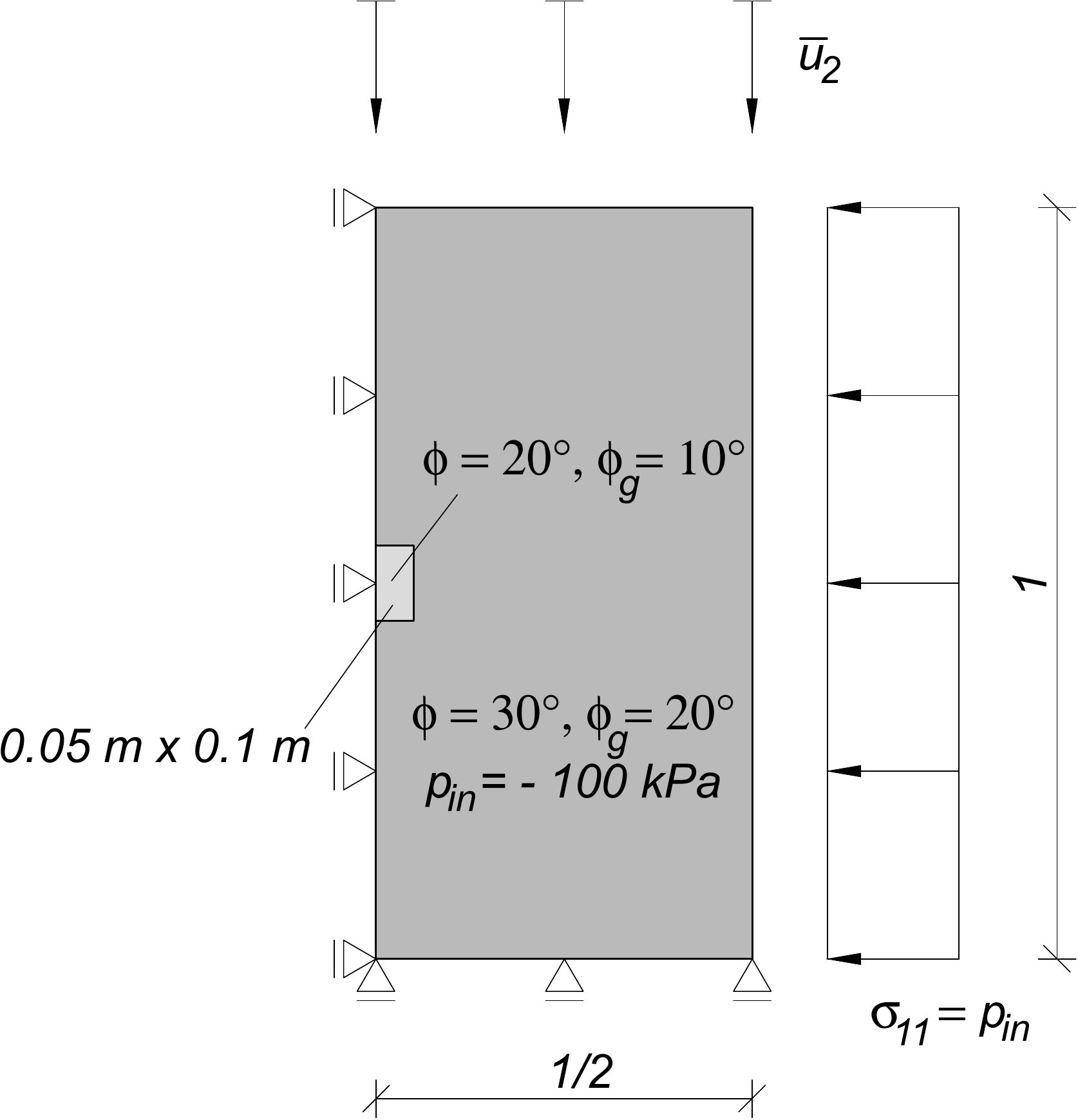}
\caption{Geometry of the domain used for analysing a  biaxial  compression test.  A defect inclusion is located in the middle of the specimen.}.
\label{Fig_BiaxialGeom}
\end{figure}

A first set of analyses was performed to simulate a biaxial compression test on a specimen of an elastic-perfect plastic, non-associated,  Mohr-Coulomb material,  characterized by a defect inclusion in its middle. 
Due to the symmetry of the problem,  only half of the specimen was analysed, and its geometry is shown in Fig. \ref{Fig_BiaxialGeom} together with the  boundary conditions.  
A vertical plane strain compression was applied,  the initial stress being spherical. 
No boundary conditions were imposed to the Cosserat rotational degrees of freedom,  the micro-volumes being   hence allowed   to rotate freely.

The analyses were conducted with  three progressively finer  meshes and both on the Cauchy and the Cosserat continuum, the same set of common parameters being used which is shown  in Table \ref{Tab_Biaxial}.
However in order to prevent a very early crush of the analyses involving the Cauchy continuum, a moderate non-associativity was used,  characterized by an angle $\phi_g$ of the plastic potential $33\%$ lower than that used to define the yield surface $\phi$.  
The parameters of the  defect inclusion  only differ for the two angles $\phi$ and $\phi_g$ and are shown in Fig. \ref{Fig_BiaxialGeom}.
A rounded Mohr-Coulomb model was used, circumscribed to the original one which is otherwise characterized by discontinuities in the deviatoric plane. However, the rounding parameter $\beta$ of both the yield and the plastic potential surfaces was very close to unity,  so that the  surfaces practically coincide with those of the original criterion. 

\begin{figure}[tb]
\centering
\includegraphics[height=0.5\textwidth]{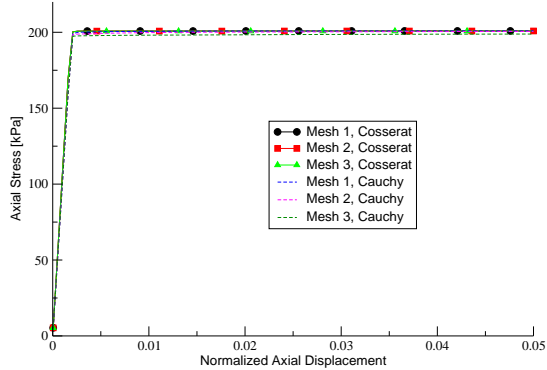}
\caption{Biaxial compression test analyses of an elastic-perfect plastic,  non associated Mohr-Coulomb material  with a moderate non associativity and a defect inclusion: Cauchy vs Cosserat continua.  Stress-strain curves.}
\label{Fig_BiaxialStressStrain}
\end{figure}

The stress-strain curves in the axial direction, for the six analyses are shown in Fig. \ref{Fig_BiaxialStressStrain}.  
No structural softening was observed in the behaviour of the Cauchy material. 
The stress-strain behaviour of both the Cauchy and the Cosserat continuum is very similar, with the six curves very close to one another. 
However whilst for the Cauchy material a slight decrees of the failure load is observed as the mesh is refined,  the three curves of the Cosserat model are perfectly coincident.

\begin{figure}[ht]

\begin{subfigure}{.3\textwidth}
  \centering
  \includegraphics[height=1.5\linewidth]{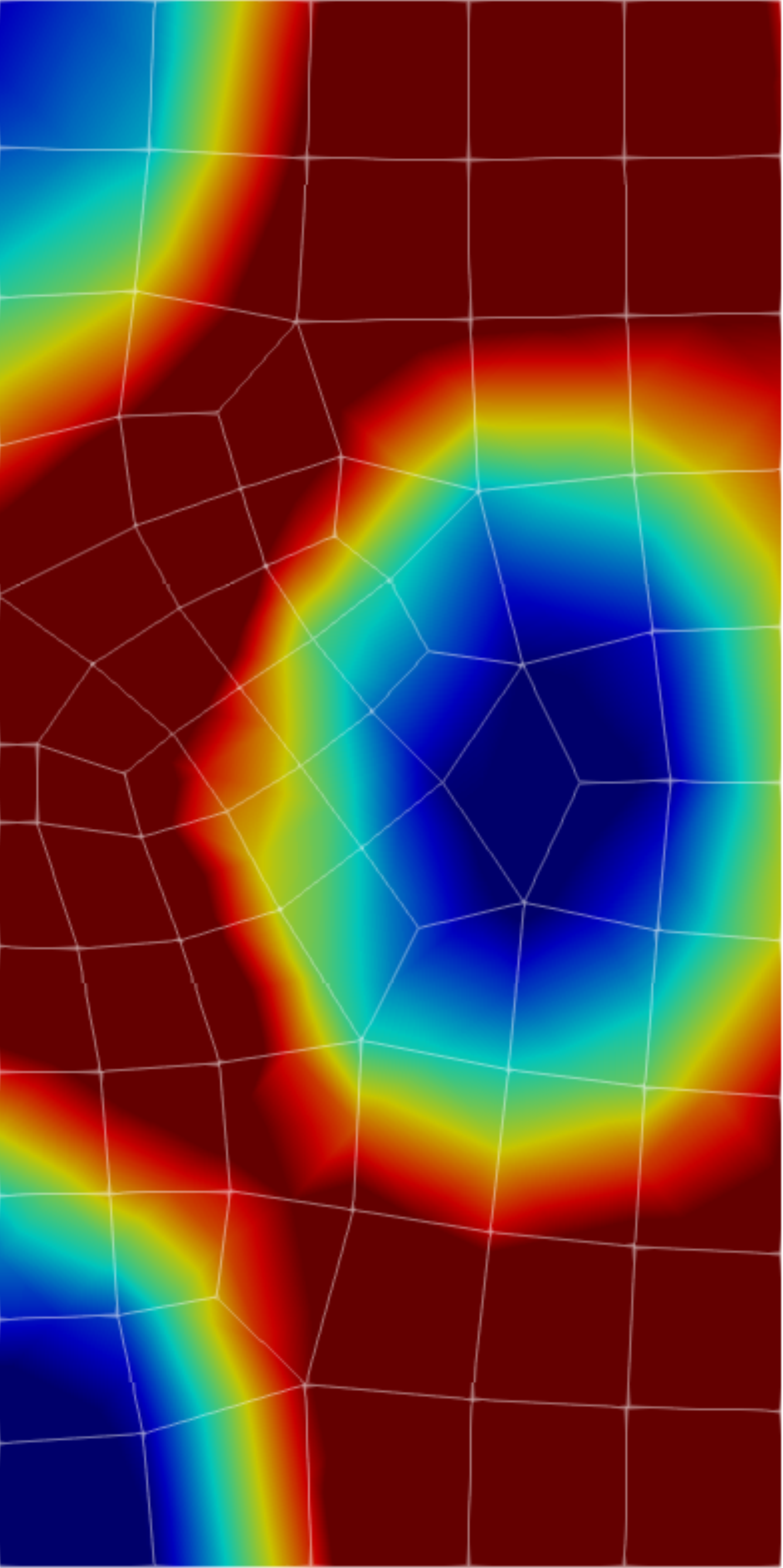}  
  \caption{Cauchy: mesh 1}
  \label{fig:sub-first}
\end{subfigure}
\begin{subfigure}{.3\textwidth}
  \centering
  \includegraphics[height=1.5\linewidth]{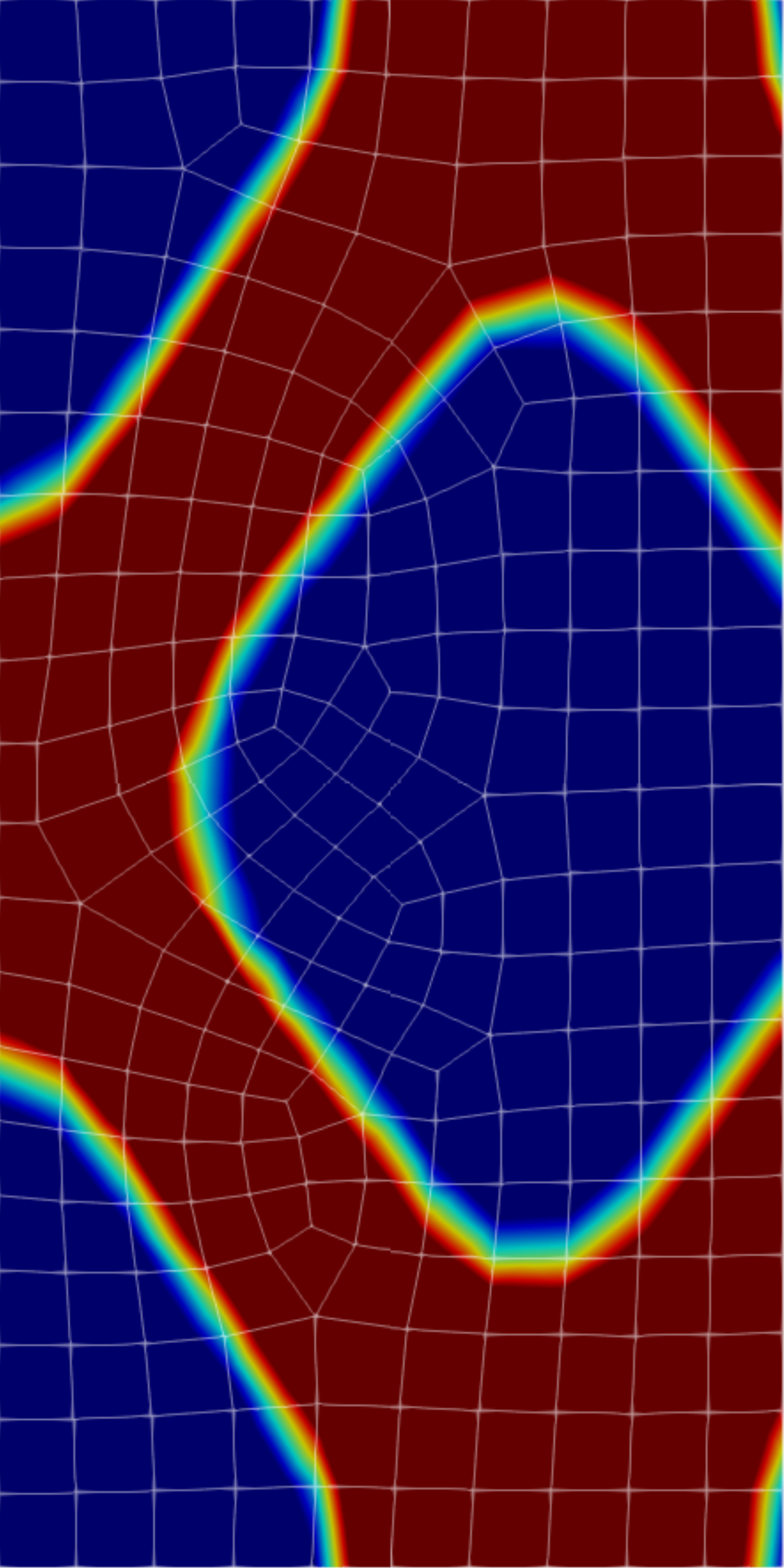}  
  \caption{Cauchy: mesh 2}
  \label{fig:sub-second}
\end{subfigure}
\begin{subfigure}{.3\textwidth}
  \centering
  \includegraphics[height=1.5\linewidth]{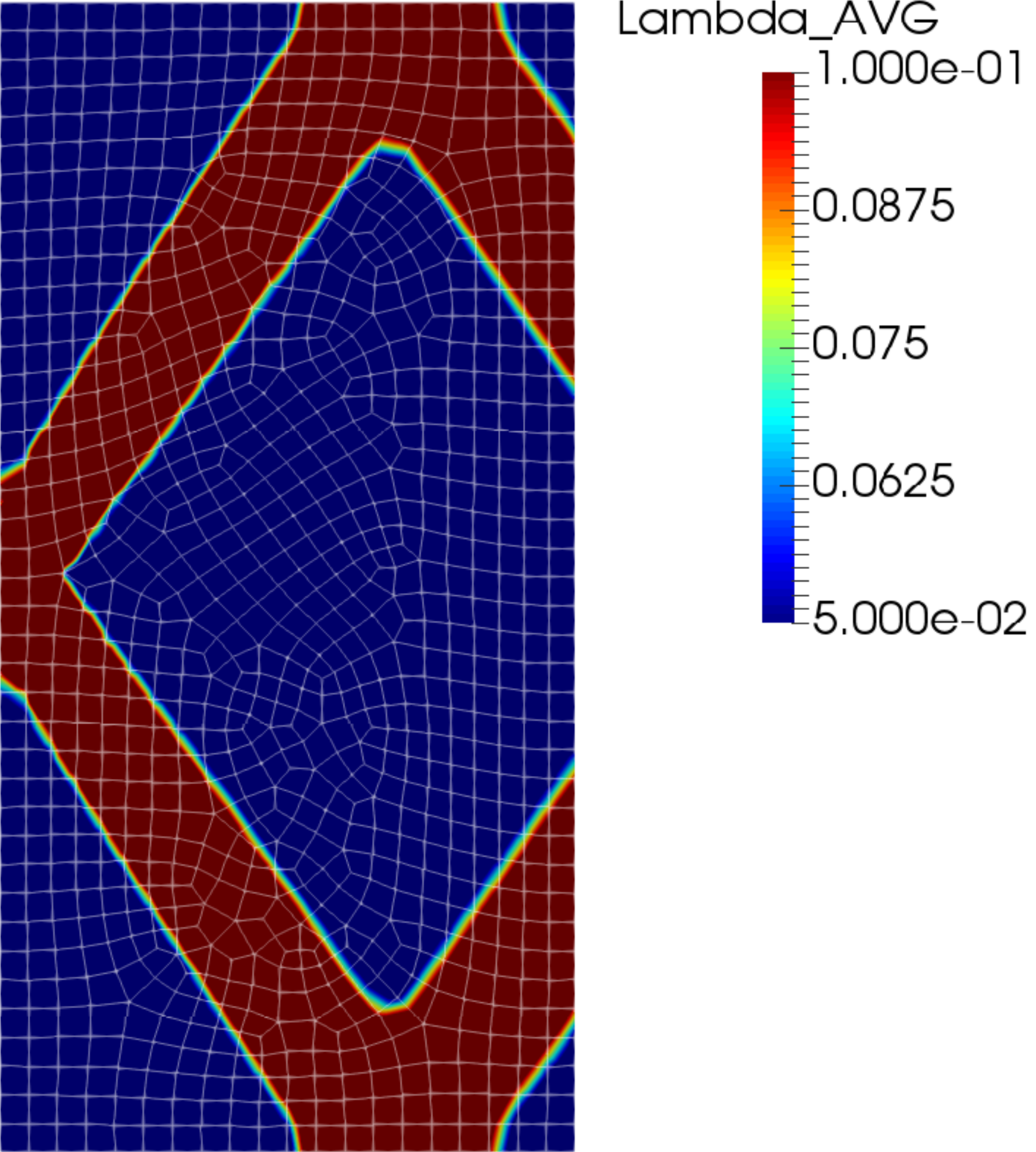}  
  \caption{Cauchy: mesh 3}
  \label{fig:sub-second}
\end{subfigure}
\newline
\begin{subfigure}{.3\textwidth}
  \centering
  \includegraphics[height=1.5\linewidth]{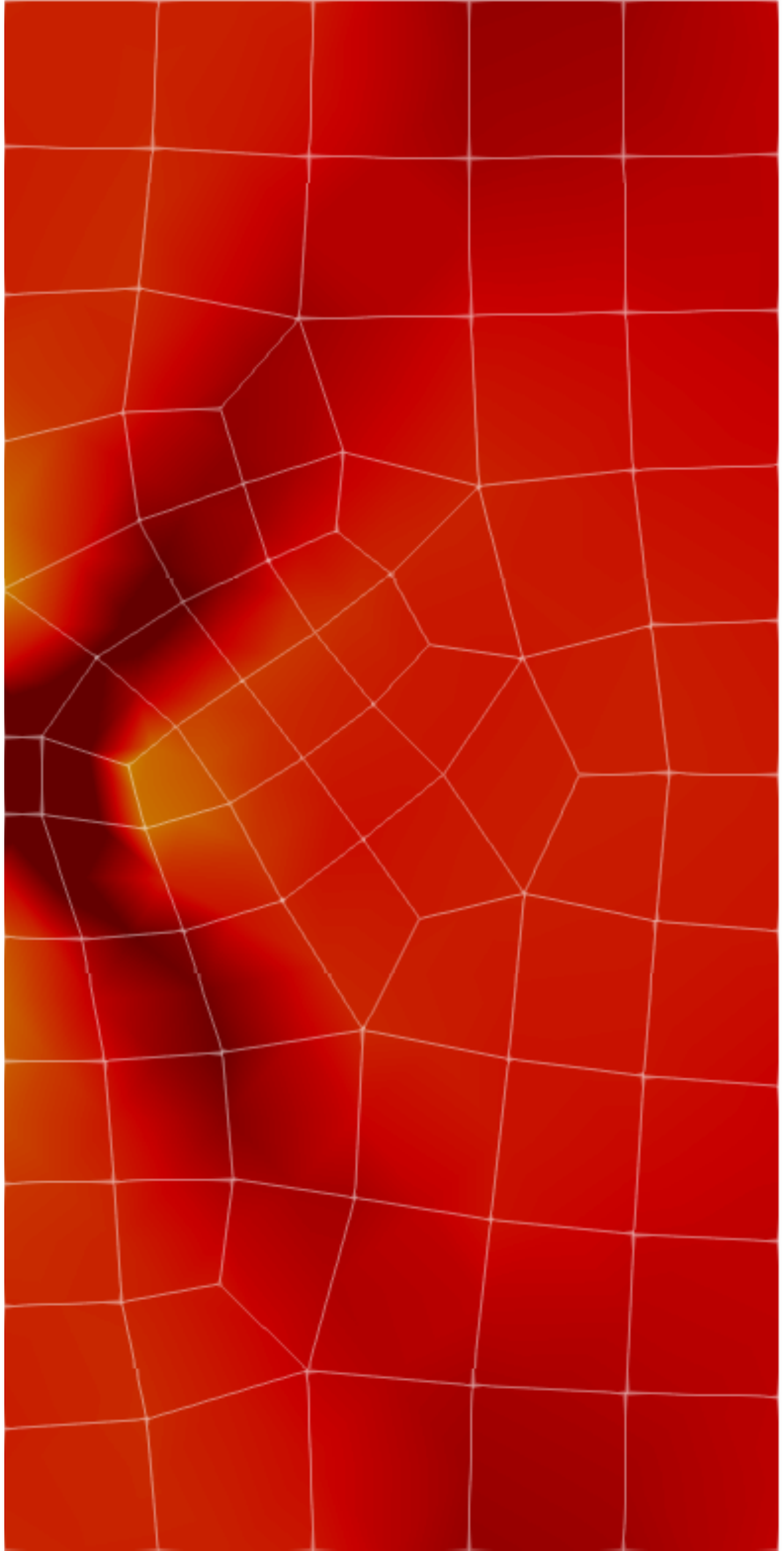}  
  \caption{Cosserat: mesh 1}
  \label{fig:sub-first}
\end{subfigure}
\begin{subfigure}{.3\textwidth}
  \centering
  \includegraphics[height=1.5\linewidth]{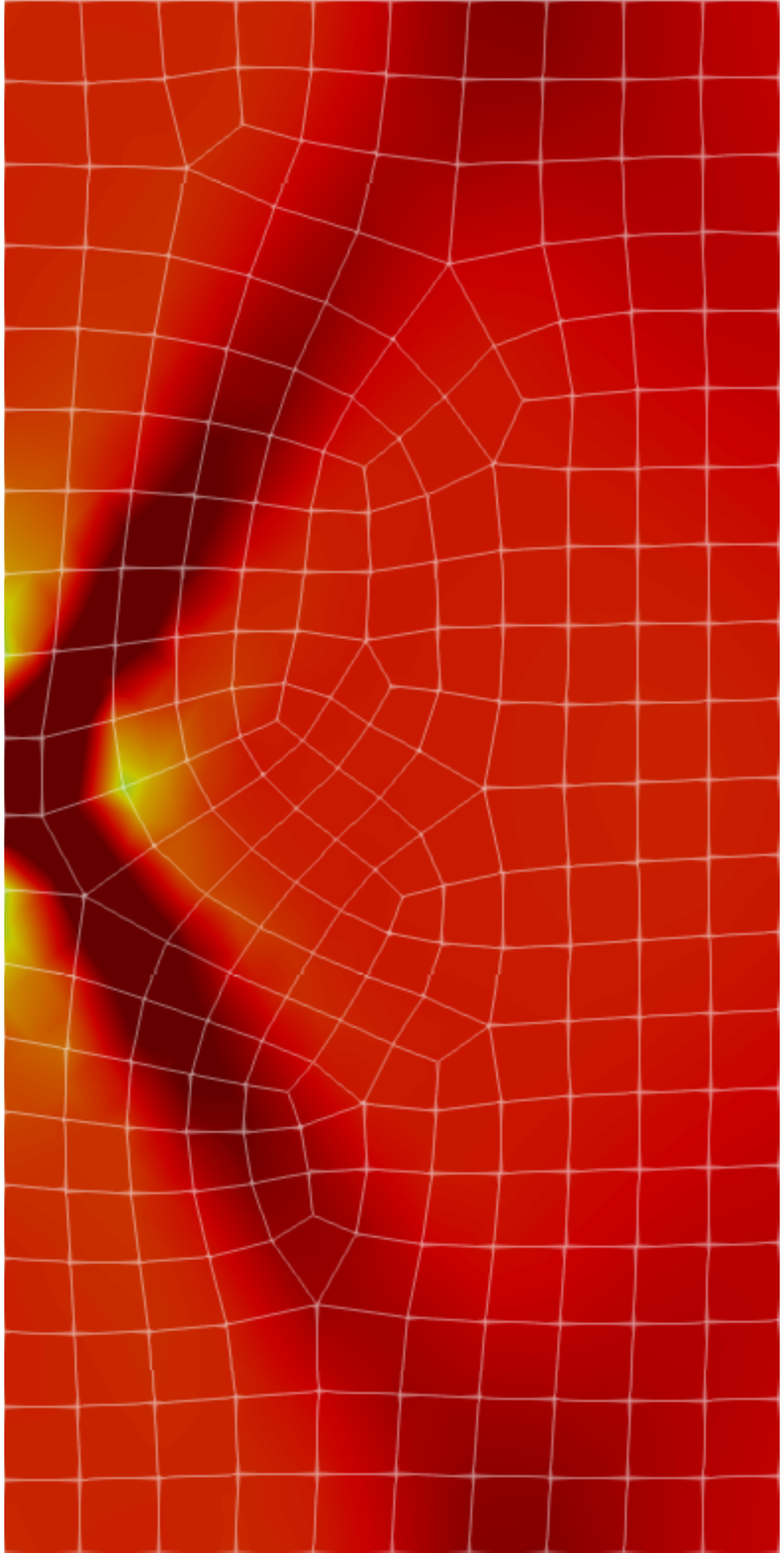}  
  \caption{Cossserat: mesh 2}
  \label{fig:sub-second}
\end{subfigure}
\begin{subfigure}{.3\textwidth}
  \centering
  \includegraphics[height=1.5\linewidth]{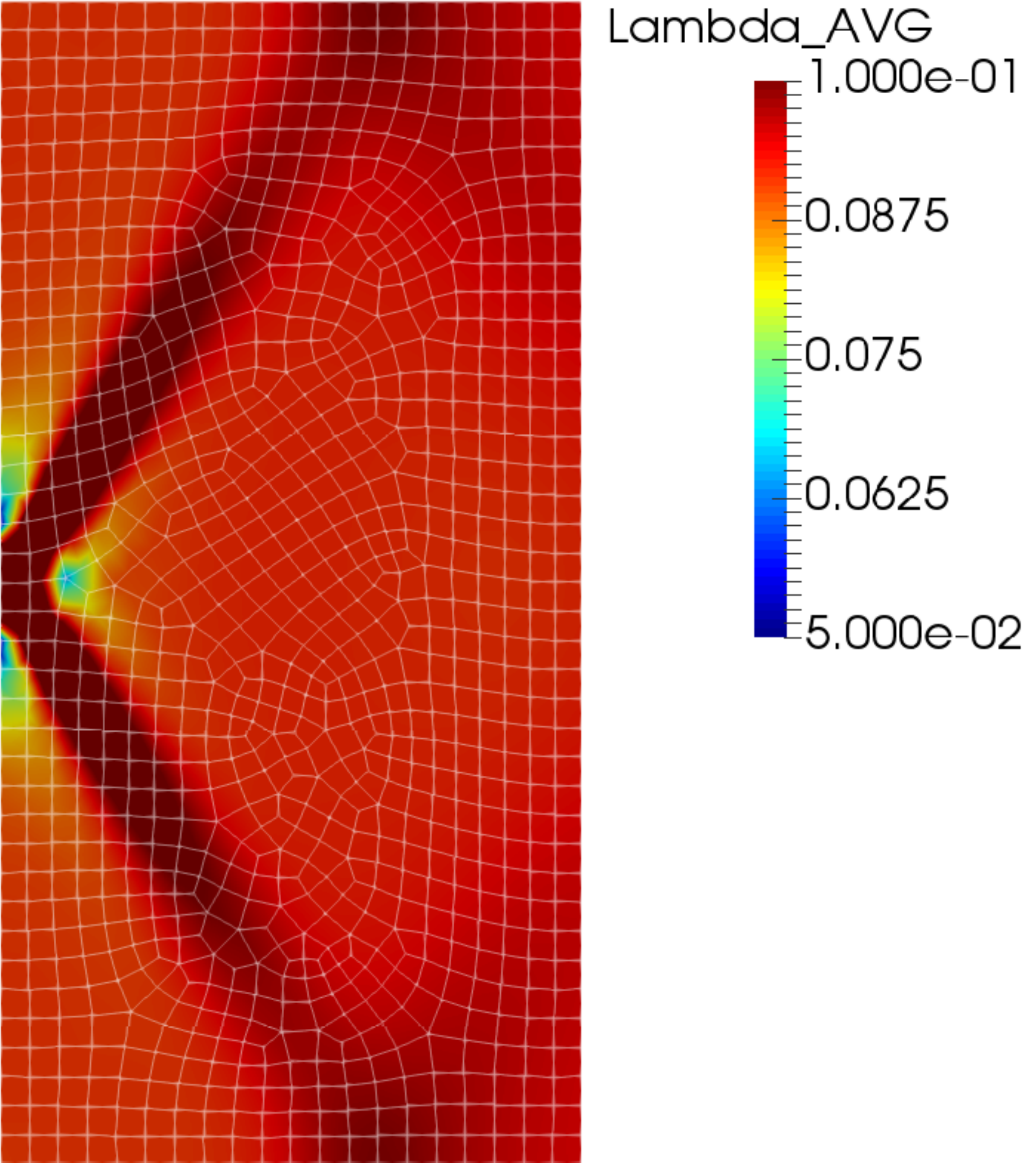}  
  \caption{Cosserat: mesh 3}
  \label{fig:sub-second}
\end{subfigure}

\caption{Biaxial compression test analyses of an elastic-perfect plastic,  non associated Mohr-Coulomb material: Cauchy vs Cosserat continua.  Contour plots of accumulated  average plastic plastic multiplier  $\lambda$. }
\label{Fig_BiaxialShearBands}
\end{figure}

Figure \ref{Fig_BiaxialShearBands} shows contours of the  plastic multiplier $\lambda$ for the six analyses.  The capability of the Cosserat continuum to eliminate mesh sensitivity is clearly apparent.

\begin{table}[t]
        \begin{tabular}{ c c c c c}
                \hline 
                \hline
                Elastic & Constit.& Yield & Plastic & Hard.\\ constants & model &function & potential & law
                \\ 
                \hline 
                $G = 55000$ kPa&      &                &                                   &  $c_i=0$\\
                $K = 33333$ kPa& Outer-Mohr & $\phi=30^\circ$&\multirow{2}{*}{$\phi_g=20^\circ$} &$c_f=0$ \\ 
                $G_c = 5000$ kPa& Coulomb& $\beta_f=0.9999$& &$a_\lambda=0$\\
                $B=B_c = 5000$ kN& & & &\\
                \hline 
                \hline
        \end{tabular}
        \caption{Parameters for biaxial test simulation}
        \label{Tab_Biaxial}
\end{table}

\section{Analysis of  a shallow  strip footing}

A shallow and rigid footing has been analysed in plane strain conditions using either the Tresca  or the Mohr-Coulomb, Matsuoka-Nakai and Lade-Duncan models.  Whilst the former is typically used to simulate undrained pore water conditions, the latter set of models is adopted for the drained regime.

A $2$m wide footing was subjected to vertical and  centred loading conditions.  Due to the symmetry of the problem only  half of the footing,  i.e.  $1m$ wide,  was modeled together with      a conservatively large   $50$m $\times50$m  soil domain   to avoid   any border effect.  

\begin{figure}
\includegraphics[width=\textwidth]{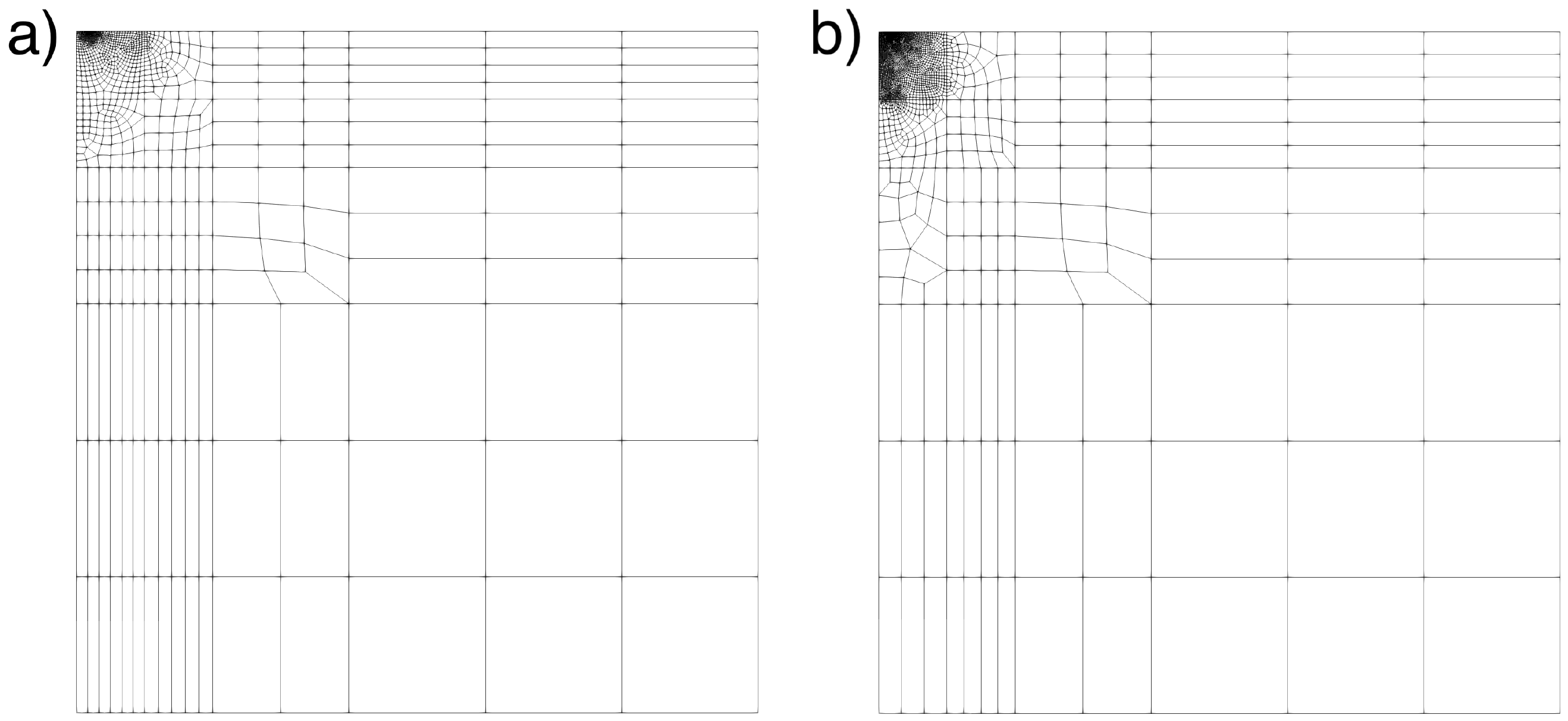}
\caption{Meshes employed to simulate the strip footing under plane strain conditions: coarse mesh (a) and fine mesh (b)}
\label{fig:meshes}
\end{figure}

The domain was discretized adopting a coarse and a fine mesh,   characterized by   3639 nodes and 1164  elements and 10009 nodes and 3272 elements, respectively (Fig. \ref{fig:meshes}).  All elements were  quadratic with eight-nodes each.
The stiffness matrix and the stresses have been fully integrated for the Cosserat continuum  whilst a selective integration was used for the classical Cauchy medium. 
All analyses  have been performed under displacement control and 
model parameters  are shown  in Table \ref{Tab_FottingParameters}. It should be noted that the rounding parameter $\beta$ for the yield and the plastic potential surfaces was set for all analyses close enough to unity, so that the rounded deviatoric sections  were very close to those of the original Tresca and Mohr-Coulomb surfaces. 
Moreover,  since the analyses were conducted in plane strain conditions,  the wryness tensor $\tenss{\chi}$ is deviatoric.  This is the consequence of the restrained variation along the out of plane direction of the only available component of the rotation vectors (which is  along the same direction).  Hence the  Cosserat elastic parameter $T$  plays no role in the analyses.

\begin{table}[t]
        \begin{tabular}{c c}
                \hline 
                \hline
                Common params, & \multirow{2}{*}{$G = 4166.7$ MPa, $K = 5555.6$ MPa}\\
                Cauchy & \\
                \hline
                Common params, & $G = 4166.7$ MPa, $K = 5555.6$ MPa,\\
                Cosserat & $G_c = 250$ MPa, $B=B_c = 250$ MN\\          
                \hline
                Strip footing & $\beta_f=0.999$, $c_i=490$ kPa, $c_f=0$ kPa,  \\
                Tresca soil & $a_\lambda=10$ (softening), $a_\lambda=0$ (perf.plast.),  \\
                \hline
                $N_c - N_\gamma$ & $\phi=25^\circ$, $c_i=464.49$ kPa, $\gamma=18$ kN/m$^3$, \\
                problem & $c_f=0$ kPa, $a_\lambda=0$ ($\phi_g=10^\circ$ non-associated flow)  \\
                \hline
                \multirow{3}{*} {$N_\gamma$ problem} & $\phi=25^\circ$, $c_i=0$ kPa, $\gamma=18$ kN/m$^3$, \\
                 & $c_f=0$ kPa, $a_\lambda=0$ ($\phi_g=0.5^\circ$ non-associated flow),  \\
                     & ($\beta_f=0.9999$ for Outer-Mohr-Coulomb)\\
                \hline
                \hline
        \end{tabular}
\caption{Parameters for shallow strip footing analyses.}
\label{Tab_FottingParameters}
\end{table}

\subsection{Tresca soil}
Analyses were conducted for both the Cauchy and the Cosserat continua adopting the Tresca criterion to describe both the yield and the plastic potential surfaces,  so that an associated flow occurred.  
A perfect plastic or a softening behaviour was considered, the latter being defined  by a high rate  exponential rule as  described in \textit{Part I} of this paper. 
No lateral surcharge and a frictionless interface at the soil-footing contact were assumed.

\begin{figure}[tb]
\centering
\includegraphics[width=0.8\textwidth]{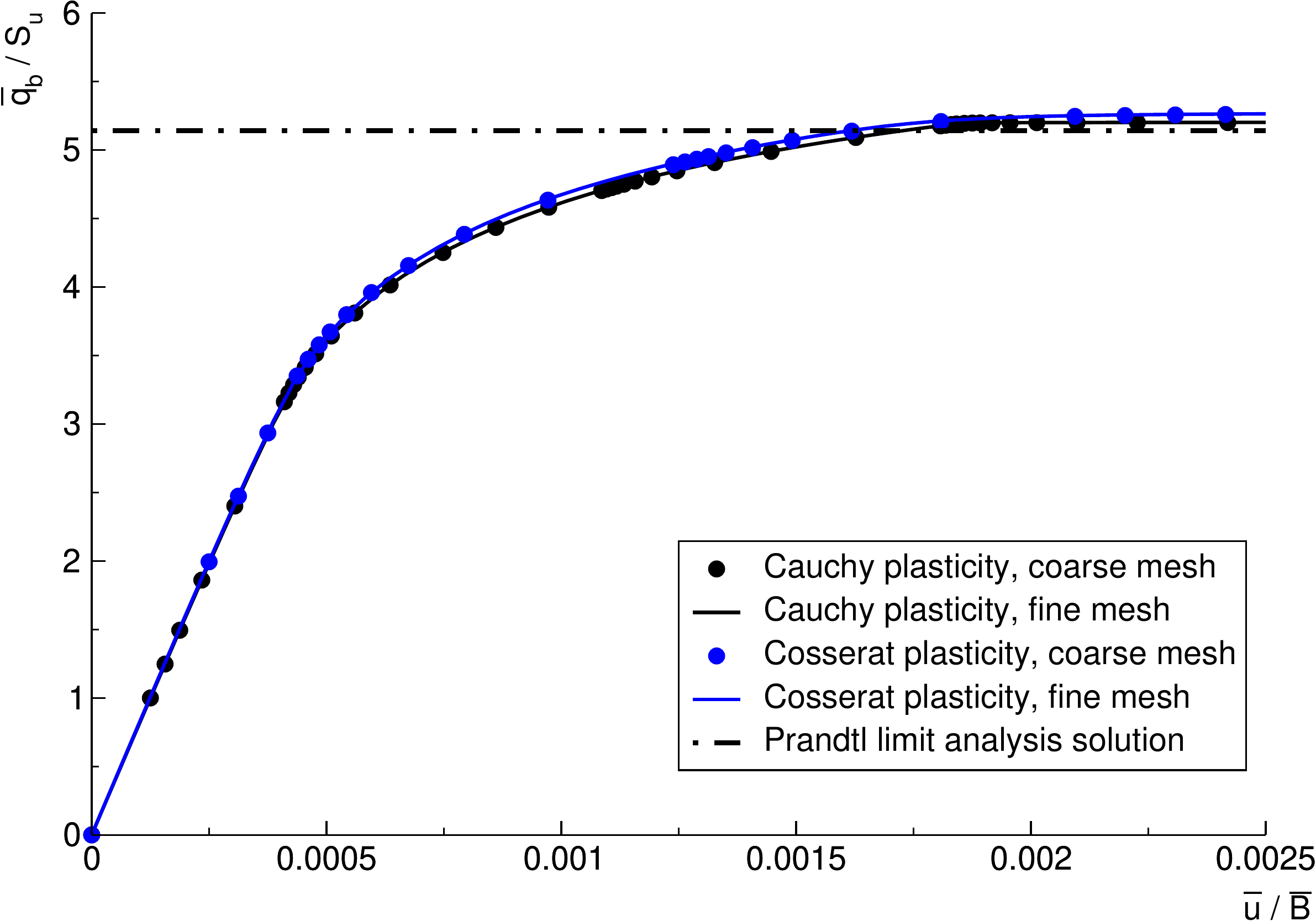}
\caption{Strip footing, Tresca soil and perfect plasticity, load-displacement curve}
\label{Fig_TrescaPerfectPlasticity}
\end{figure}

Fig. \ref{Fig_TrescaPerfectPlasticity}  shows the load-displacement curves for analyses conducted in associated perfect plasticity for the Cauchy and Cosserat continua. 
 The bering pressures $\bar{q}_b$ and the vertical displacements $\bar{u}$ have been normalized with respect to the undrained shear strength $S_u$ and the breath $\bar{B}=2$m  of the foundation, respectively.
The dash-dotted line in the  diagram  shows the theoretical bearing capacity factor
\begin{equation}
N_c= \frac{\bar{q}_{f}}{S_u}= 2+\pi
\nonumber
\end{equation}
of the Prandtl's  solution  \cite{Prandtl1920} valid for associated perfect plasticity in the absence of a lateral surcharge.

\begin{figure}[tb]
\centering
\includegraphics[width=0.8\textwidth]{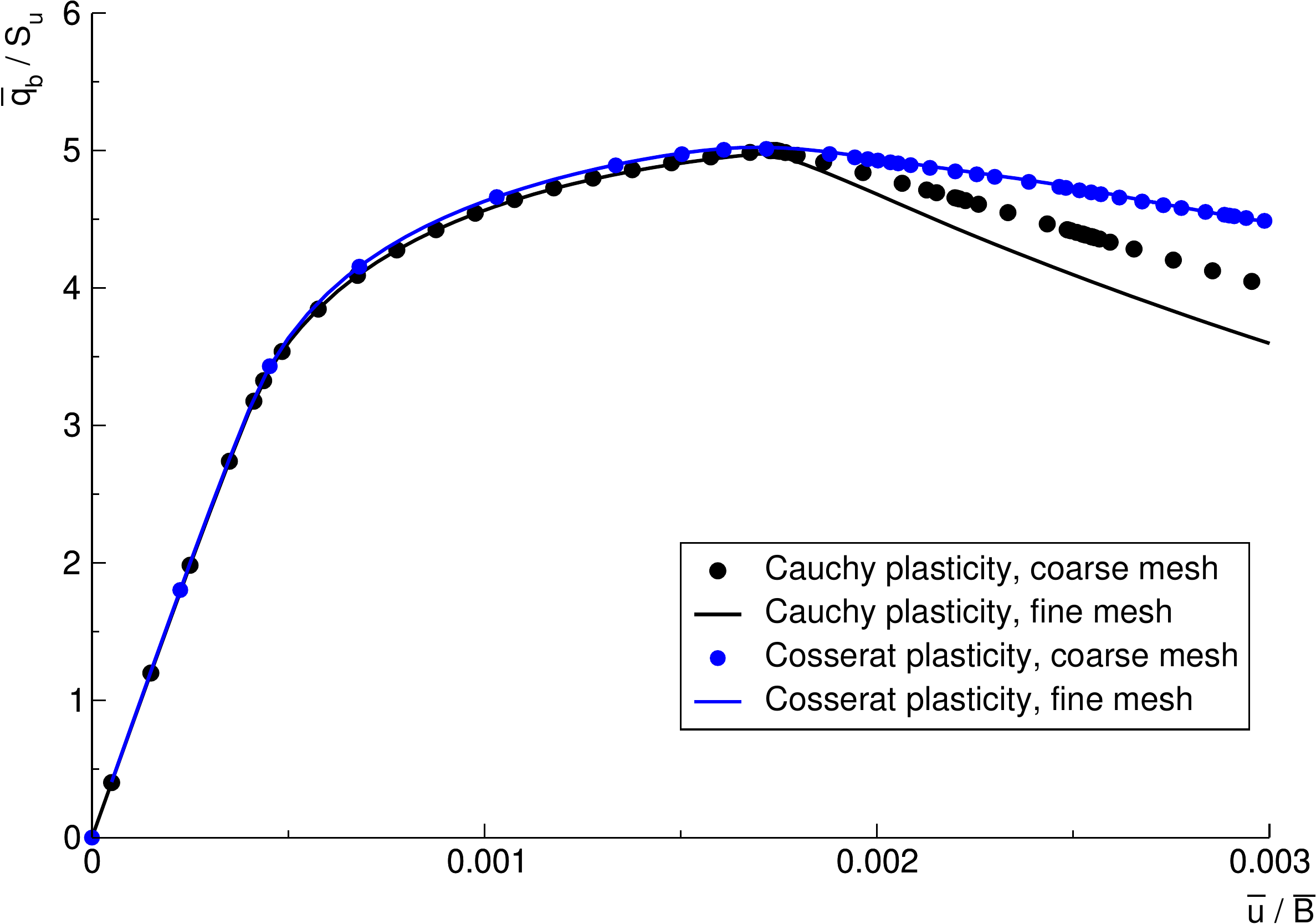}
\caption{Strip footing, Tresca soil and exponential softening law, load-displacement curve}
\label{Fig_TrescaSoftening}
\end{figure}

Fig. \ref{Fig_TrescaSoftening} shows the  load-displacement curves of the analyses conducted in associated plasticity and exponential softening,  for   the Cauchy and Cosserat continua using  both coarse and fine meshes.

%
%

As expected,  Fig. \ref{Fig_TrescaPerfectPlasticity} clearly indicates that there is no appreciable mesh dependency when analyses are conducted adopting  perfect and associated  plasticity.  The result is confirmed for both the Cauchy and the Cosserat continua.

However  when softening is included in the Cauchy continuum  different load-displacement curves are obtained for the coarse and fine meshes,  showing as expected the existence of strong mesh dependencies (Fig. \ref{Fig_TrescaSoftening}).
No such effect exists for the Cosserat continuum,  since the curves of the two analyses conducted with the coarse and the fine meshes perfectly coincide,  confirming the effectiveness of such a continuum to provide a mesh independent solution.

In the case of perfect associated plasticity the results of the analyses provide a bearing capacity factor very close to $N_c=2+\pi$  of  the theoretical Prandtl solution.  As expected both the Cauchy and the Cosserat continua slightly overestimate that value by  $1.1\%$ and $2.4\%$ respectively (Fig. \ref{Fig_TrescaPerfectPlasticity}).

When an exponential softening is included (Fig. \ref{Fig_TrescaSoftening}) all four analyses for the coarse and fine meshes and for the Cauchy and the Cosserat are very close to one another until a peak is attained. Thereafter they diverge,  the latter continuum exhibiting the least reduction in  bearing capacity factor.  
The reduction in bearing capacity is presumably associated to a progressive failure effect caused by the exponential reduction in the undrained strength.

\subsection{Mohr-Coulomb,  Matsuoka-Nakai and Lade-Duncan soils}
The so-called \textit{pure} $N_\gamma$ problem was investigated.  The yield and the plastic potential surfaces were defined by the Mohr-Coulomb,  Matsuoka-Nakai and Lade-Duncan criteria, both in associated and non-associated plasticity, the latter being achieved by using a different set of parameters for the plastic potential. No hardening/softening behaviour was included.
The horizontal displacement of the footing nodes was restrained to simulate a perfect rough interface. 
 
\begin{figure}[tb]
\centering
\includegraphics[width=0.8\textwidth]{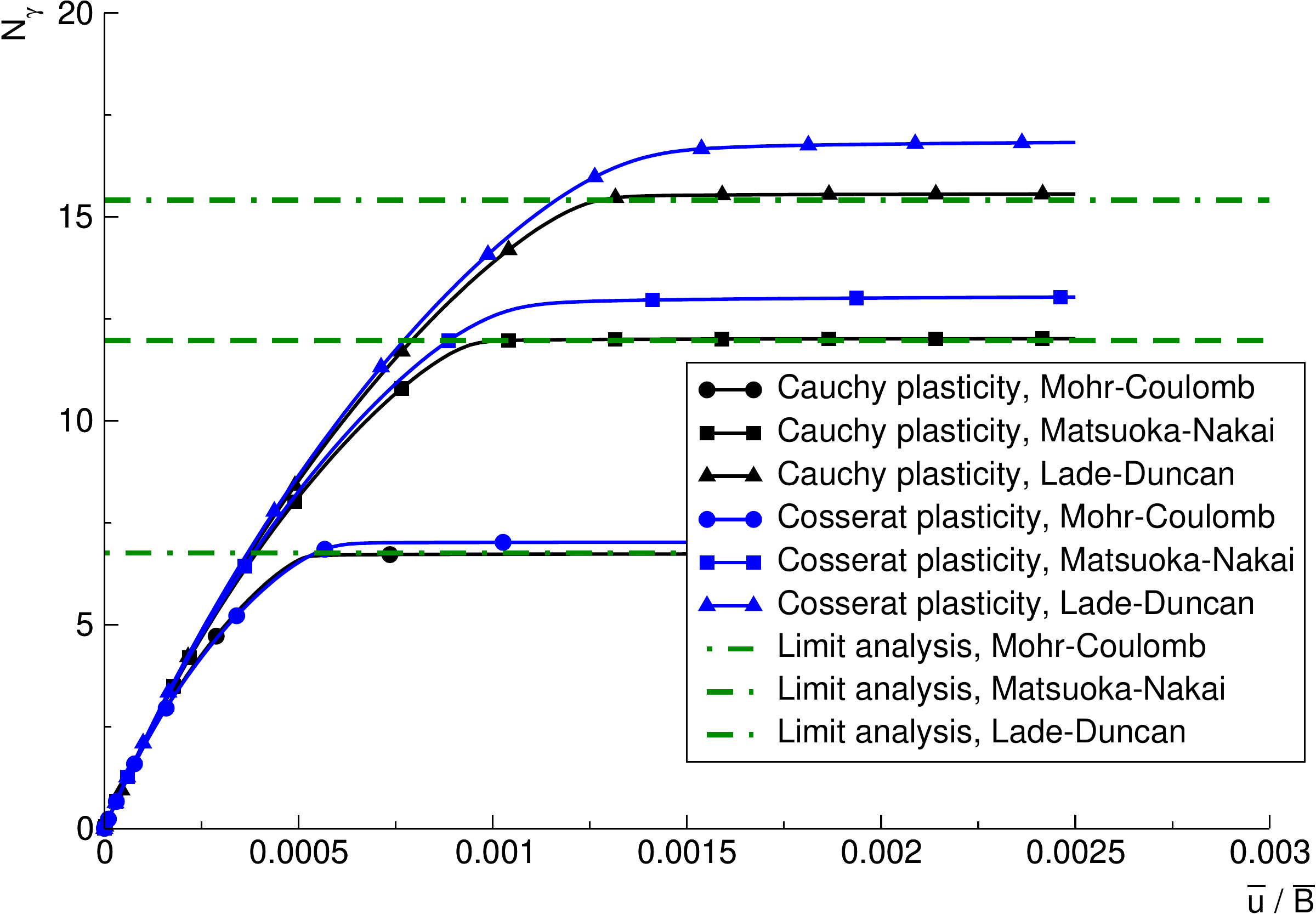}
\caption{Strip footing on Mohr-Coulomb,  Matsuoka-Nakai and Lade-Duncan soils in associated plasticity.}
\label{Fig_NgammaAss}
\end{figure}

Figure \ref{Fig_NgammaAss} shows a comparison between analyses conducted on the Cauchy and the Cosserat continua,  in associated perfect plasticity.  Theoretical values of $N_\gamma$ for the three constitutive models are also shown for comparison. These were obtained by means of the limit analyses using the \textit{ABC} software \cite{Martin2004} in combination with  the findings of Lagioia and Panteghini \cite{Lagioia2017}.
Results indicate that the bearing capacity in  Cauchy perfect plasticity coincides with that obtained from the limit analyses for all three failure criteria, whilst for the Cosserat continuum  consistently larger values for all models were observed.  However for the Mohr-Coulomb criterion the difference between the Cauchy and the Cosserat bearing capacity was not significant.

\begin{figure}[tb]
\centering
\includegraphics[width=0.8\textwidth]{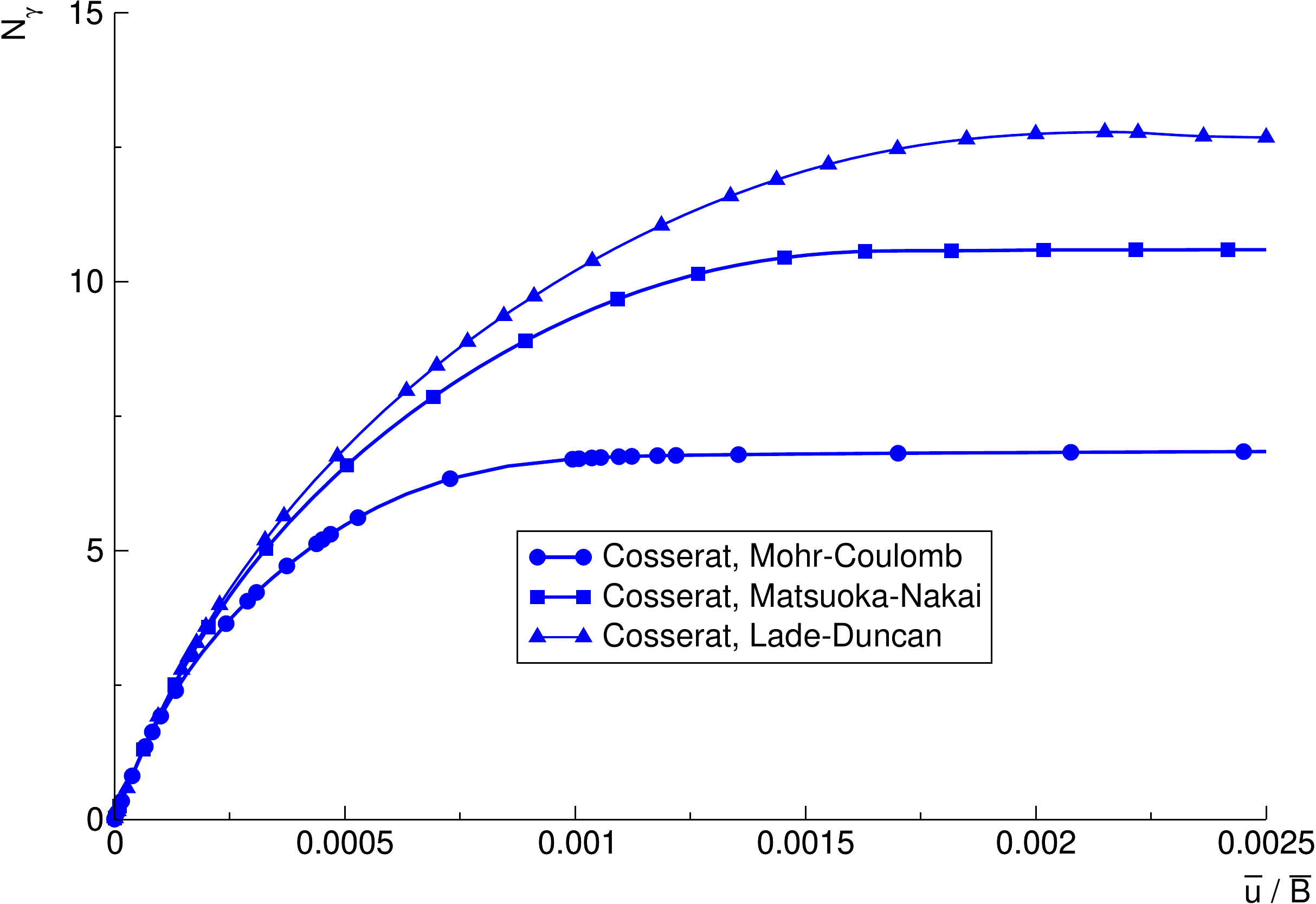}
\caption{Strip footing on Mohr-Coulomb,   Matsuoka-Nakai and Lade-Duncans Cosserat  soils with a  limit non-associated plasticity.}
\label{Fig_NgammaNonAss}
\end{figure}

The  same set of analyses was also conducted with an extreme non-associated plastic flow,  by setting the angle $\phi_g$ of the plastic potential virtually to zero  ($\phi_g=0.5^\circ$). 
As expeteced, the analyses on the Cauchy continuum crushed straight away, without performing even the first load step. On the contrary no problems were encountered with the Cosserat material. Moreover all analyses were concluded in just a few minutes on a portable computer. The $N_\gamma$-displacement curves are shown in Fig. \ref{Fig_NgammaNonAss}.

A final set of analyses was conducted on the Matsuoka-Nakai soil, with the aim of comparing associated and non-associated plasticity  for both the Cauchy and the Cosserat continuum.  In order to manage to avoid crushing of the analyses involving the former medium,  a moderate non-associativity was used,  setting the difference between  the angles $\phi$ and $\phi_g$  to $15^\circ$.  As this remedy alone was not enough,   an effective cohesion  was also introduced in the yield surface (see Table \ref{Tab_FottingParameters}).

\begin{figure}[tb]
\centering
\includegraphics[width=0.8\textwidth]{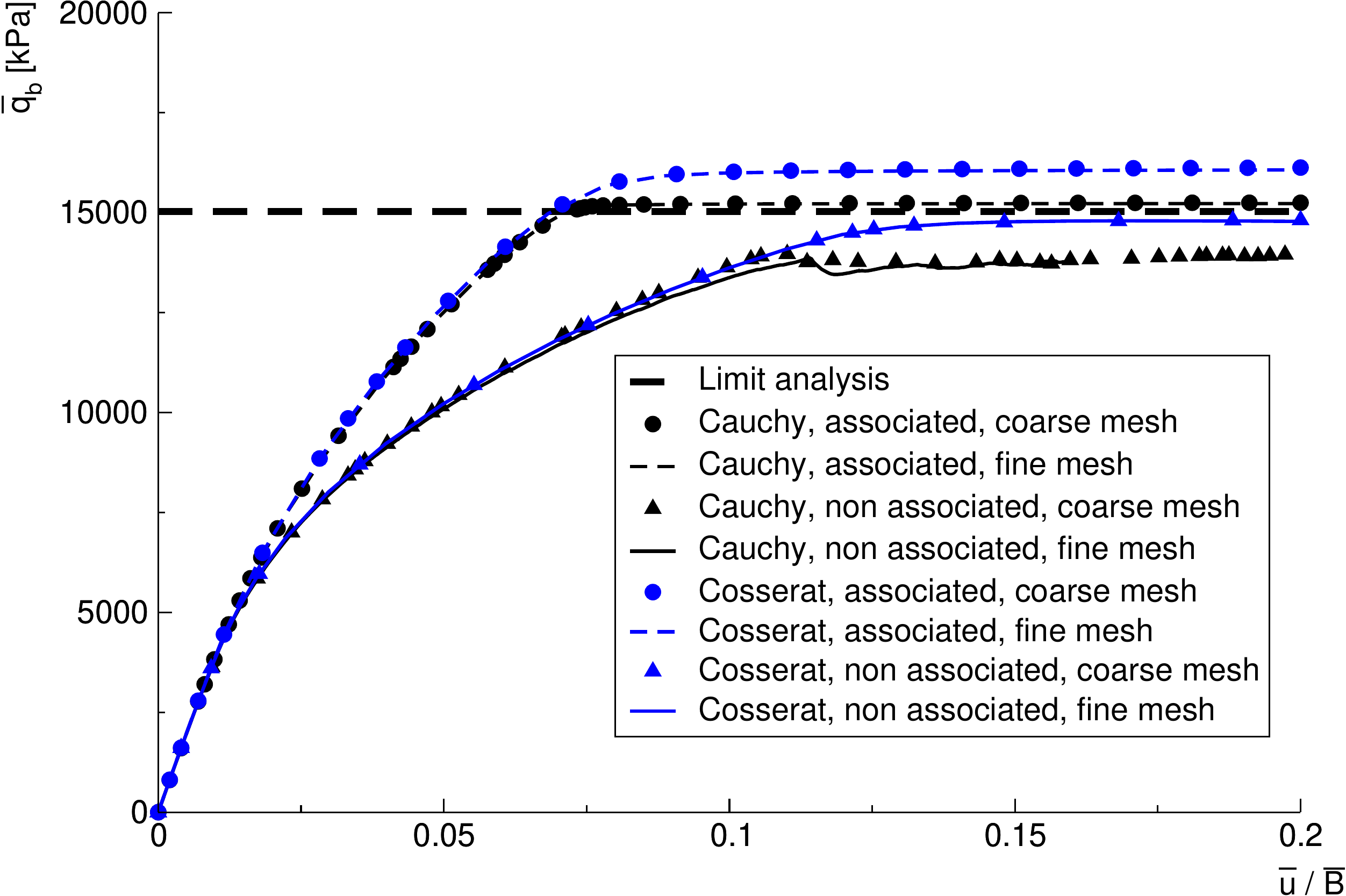}
\caption{Strip footing on Matsuoka-Nakai soil with moderate non-associativity and cohesion.}
\label{Fig_fondaz_mn_coes}
\end{figure}

Figure \ref{Fig_fondaz_mn_coes} shows that no mesh sensitivity occurs with the Cosserat continuum  since the analyses of the coarse and the fine meshes yielded identical results  both in associated and non-associated plasticity.
As expected,  the same conclusion can be drawn for the  Cauchy medium  in associated plasticity.  However mesh sensitivity and structural softening emerged when the flow was non associated.

\section{Conclusions}
A Finite Element procedure for the integration of the 
 elasto-plastic constitutive model described in \textit{part I} of this paper for the Cosserat continuum has been   implemented using a full backward Euler predictor/corrector scheme into a proprietary Finite Element program.
The integration algorithm is an extension to the Cosserat medium of that proposed by Panteghini and Lagioia \cite{PL2018} for the classical Cauchy/Maxwell continuum.

At variance with other approaches available in the literature, which adopt as principal variables either the six components of the strain tensor or its principal values,  the Panteghini and Lagioia scheme is based on invariants and results in a number of important benefits. 
First of all the integration requires  the solution of a single equation in a single unknown, which brings considerable speed improvements as compared to the system of seven by seven and four by four equations and unknowns necessary for  the other two approaches.
Moreover no difficulties are encountered when two and/or three of the principal stresses coincide. 
Finally the use of invariants anables the adoption of the technique proposed by De Souza Neto et al. \cite{bibbia} to deal with the singularity at the apex of the yield and plastic potential surfaces.

The drawback of the adoption of an implicit integration with the  Cosserat continuum is the large amount of calculations required,  since nine consistent tangent operators needs to be derived.
However since the constitutive model proposed in \textit{part I} of this paper was formulated to allow the use of multiple classical yield criteria,  which are relevant in both research and engineering practice,  and since these criteria are highly non-linear,  the choice of an implicit algorithm is recommended.  This is in fact the only one which ensures unconditional stability with highly non-linear surfaces. Moreover the use of consistent tangent operators results in a quadratic convergence of the structural Newton loop. 

The reward of the adoption of the proposed integration scheme is considerable.  Analyses with a large number of elements and nodes  were completed within a couple of minutes and no numerical difficulties were ever encountered,  despite the fact that a very demanding boundary value problem was investigated. 

A first set of analyses was conducted to model a biaxial compression test in plane strain conditions on a Mohr-Coulomb elastic-perfect plastic material with non-associated flow.  A defect inclusion was included in the middle of the specimen.  The comparison of the Cauchy and Cosserat continuua, clearly confirm the capability of the latter to remove mesh sensitivity.  

A plane strain shallow footing problem has been considered next,  the soil behaviour being reproduced either with a Tresca or with a Mohr-Coulomb,  Matsuoka-Nakai and Lade-Duncan models.  
The analyses conducted with the latter model indicate that the Cauchy and the Cosserat continuum exhibit the same load-displacement curves, with no mesh sensitivity,  when a perfect plastic and associated behaviour is considered. 
However the introduction of a very non linear  softening  results in mesh sensitivity for the Cauchy medium, whilst the such effect emerges for the Cosserat.

The set of analyses performed with the Mohr-Coulomb, Matsuoka-Nakai and Lade-Duncan materials shows that the Cosserat continuum is effective in removing mesh sensitivity,  structural softening  and early crush of the analyses when non associated flow is considered.  

%

\bibliographystyle{plain}
\bibliography{short,paper}

\begin{appendices}

\section{Finite Elements formulation}
The proposed Finite Elements (FE) formulation is based on a spatial discretization of the fields $\tensf{u}$ and $\tensf{\theta}$, whose nodal values are contained in the vector of the generalized displacements $\hat{\tensf{u}}_g$. 
In the general  three-dimensional case   it is a vector of $6\hat n$ components, $\hat n$ being the total number of nodes of the element,  which reduces to   $3 \hat n$ for planar problems.
 
The spatial discretization is based on the choice of the shape functions $N^{(i)}(\xi,\eta,\zeta)$, where $i=1,...\hat n$, with $\xi=[-1,1], \eta=[-1,1], \zeta=[-1,1]$ indicate the intrinsic parent element coordinates. 
The same shape functions are adopted to discretize both the geometry (i.e., to map the real coordinates $\tensf{x}$ with respect to the intrinsic coordinates $\tensf{\xi}$) and the generalized displacements $\tensf{u}_g$. It is then assumed:
\begin{equation}
\begin{gathered}
\tensf{x}  (\xi,\eta,\zeta)= \sum_{i=1}^{\hat n} N^{(i)}(\xi,\eta,\zeta) \hat{\tensf{x}}^{(i)}\\
\tensf{u}_g (\xi,\eta,\zeta) =\begin{bmatrix} \tensf{u}(\xi,\eta,\zeta)\\ \tensf{\theta}(\xi,\eta,\zeta) \end{bmatrix}= \sum_{i=1}^{\hat n} N^{(i)}(\xi,\eta,\zeta) \hat{\tensf{u}}_g^{(i)}
\end{gathered}
\end{equation}
where $\hat{\tensf{x}}^{(i)}$ and  $\hat{\tensf{u}}_g^{(i)}$ are respectively the vectors of the coordinates and of the generalized displacements of the $i-$th node.
Once this discretization has been adopeted, the vectors $\check{\tensf{\varepsilon}}  (\xi,\eta,\zeta) $, $\check{\tensf{\omega}}  (\xi,\eta,\zeta)$ and $\check{\tensf{\chi}}  (\xi,\eta,\zeta)$, containing the components of the tensors $\tenss{\varepsilon}(\xi,\eta,\zeta)$, $\tenss{\omega}(\xi,\eta,\zeta)$ and $\tenss{\chi}(\xi,\eta,\zeta)$ respectively can be computed as a function of the nodal values of the generalized displacaments $\hat{\tensf{u}}_g$ as:
\begin{equation}
\begin{gathered}
\check{\tensf{\varepsilon}} (\xi,\eta,\zeta) = \bvec B (\xi,\eta,\zeta) \hat{\tensf{u}}_g\\
\check{\tensf{\omega}}  (\xi,\eta,\zeta)= \bvec W  (\xi,\eta,\zeta)\hat{\tensf{u}}_g\\
\check{\tensf{\chi}} (\xi,\eta,\zeta)=\bvec M  (\xi,\eta,\zeta)\hat{\tensf{u}}_g
\end{gathered}
\label{eq:discrtens}
\end{equation}
where the matrices $\bvec B (\xi,\eta,\zeta)$, $\bvec W  (\xi,\eta,\zeta)$ and $\bvec M  (\xi,\eta,\zeta)$
are functions of the assumed shape functions and of their spatial derivatives. The complete description of these matrices for the plane strain case is reported in Appendix D.

Substituting Eqs. \eqref{eq:discrtens} into the internal energy,
\begin{equation}
\Psi= \int d\Psi=  \int \left(\tenss{\sigma}_{sym} \colon \dot{\tenss{\varepsilon}}^e + \tenss{s}_{skw} \colon \dot{\tenss{\omega}}^e + \tenss{\mu}\colon \dot{\tenss{\chi}}\right)
\label{Eq_ElasticStrainEnergyPotential_GenSymSkew}
\end{equation}
yields
\begin{equation}
W_{int}=\int_\Omega \left( \bvec B^T \check{\tensf{\sigma}}_{sym} +\bvec W^T \check{\tensf{s}}_{skw}+\bvec M^T \check{\tensf{\mu}} \right) \delta \hat{\bvec u}_g \de V
\end{equation}
where $\check{\tensf{\sigma}}_{sym}$, $\check{\tensf{s}}_{skw}$ and $\check{\tensf{\mu}}$ are the vectors containing the components of the tensors $\tenss{\sigma}_{sym}$, $\tenss{s}_{skw}$ and $\tenss{\mu}$ respectively.
To obtain the internal nodal forces $\hat{\bvec F}_{int}$ one must differentiate this last expression with respect to the variation of the nodal variables $\delta \hat{\tensf{u}}_g$
\begin{equation}
\hat{\bvec F}_{int}=\int_\Omega  \left( \bvec B^T \check{\tensf{\sigma}}_{sym} +\bvec W^T \check{\tensf{s}}_{skw}+\bvec M^T \check{\tensf{\mu}} \right) \de V
\label{eq:fint}
\end{equation}
The integration of this last equation requires the transformation between real and intrinsic coordinate system, following the standard approach  (see for instance \cite{ZTZ2005}) and is here omitted for sake of brevity.
The consistent stiffness matrix $\bvec K$ is necessary in a FE approach, and it is specifically required to assure second-order convergence of the Newton-Raphson scheme employed to solve the nonlinear algebraic system $\hat{\bvec F}_{int}=\hat{\bvec F}_{ext}$. 
It is obtained differentiating  Eq. \eqref{eq:fint} with respect to the nodal variables  $\hat{\tensf{u}}_g$, and it is equal to:
\begin{equation}
\begin{aligned}
\bvec K =\int_\Omega \left[ 
\bvec B^T \left(\der{\check{\tensf{\sigma}}_{sym}}{\check{\tensf{\varepsilon}}} \bvec B+\der{\check{\tensf{\sigma}}_{sym}}{\check{\tensf{\omega}}} \bvec W+\der{\check{\tensf{\sigma}}_{sym}}{\check{\tensf{\chi}}} \bvec M\right) \right.\\
\left.
+\bvec W^T\left(\der{\check{\tensf{s}}_{skw}}{\check{\tensf{\varepsilon}}} \bvec B+\der{\check{\tensf{s}}_{skw}}{\check{\tensf{\omega}}} \bvec W+\der{\check{\tensf{s}}_{skw}}{\check{\tensf{\chi}}} \bvec M\right)\right.\\
\left.
+\bvec M^T \left(\der{\check{\tensf{\mu}}}{\check{\tensf{\varepsilon}}} \bvec B+\der{\check{\tensf{\mu}}}{\check{\tensf{\omega}}} \bvec W+\der{\check{\tensf{\mu}}}{\check{\tensf{\chi}}} \bvec M\right)
\right]
\de V
\end{aligned}
\end{equation}
where the matrices $\derl{\check{\tensf{\sigma}}_{sym}}{\check{\tensf{\varepsilon}}}$, $\derl{\check{\tensf{\sigma}}_{sym}}{\check{\tensf{\omega}}}$, $\derl{\check{\tensf{\sigma}}_{sym}}{\check{\tensf{\chi}}}$,  $\derl{\check{\tensf{s}}_{skw}}{\check{\tensf{\varepsilon}}} $, $ \derl{\check{\tensf{s}}_{skw}}{\check{\tensf{\omega}}}$, $\derl{\check{\tensf{s}}_{skw}}{\check{\tensf{\chi}}}$, 
$\derl{\check{\tensf{\mu}}}{\check{\tensf{\varepsilon}}}$, $\derl{\check{\tensf{\mu}}}{\check{\tensf{\omega}}}$, $\derl{\check{\tensf{\mu}}}{\check{\tensf{\chi}}}$ contain the components of the fourth-order tensors 
$\derl{\tenss{\sigma}_{sym}}{\tenss{\varepsilon}}$, 
$\derl{\tenss{\sigma}_{sym}}{\tenss{\omega}}$, 
$\derl{\tenss{\sigma}_{sym}}{\tenss{\chi}}$,  
$\derl{\tenss{s}_{skw}}{\tenss{\varepsilon}}$, 
$\derl{\tenss{s}_{skw}}{\tenss{\omega}}$, 
$\derl{\tenss{s}_{swk}}{\tenss{\chi}}$,  
$\derl{\tenss{\mu}}{\tenss{\varepsilon}}$, 
$\derl{\tenss{\mu}}{\tenss{\omega}}$ and
$\derl{\tenss{\mu}}{\tenss{\chi}}$,  
 respectively. 
 These tensors are the jacobians of the stresses with respect to the strains and, in order to achieve quadratic convergence, they must be computed consistently with the stress integration algorithm.

\section{Derivatives for the general return algorithm}

\begin{equation}
\der{r}{\theta_s}=\frac{{q_s^{*}}^2 \sin\left(2\theta_s^*-2\theta_s\right)}{2 r}
\end{equation}

\begin{equation}
\der{\Delta \lambda}{\theta_s}=-\Delta \lambda\left( \frac{1}{r} \der{r}{\theta_s}+ \frac{\hat \Gamma''(\theta)}{\hat \Gamma'(\theta)}+
\frac{2}{\tan\left(2\theta_s^*-2\theta_s\right)}\right)
\end{equation}

\begin{equation}
\der{p}{\theta_s}=-\hat M K \der{\Delta \lambda}{\theta_s}           
\end{equation}

\begin{equation}
\der{q }{\theta_s}=\der{r}{\theta_s}  -3G  \left( \hat \Gamma \der{\Delta \lambda}{\theta_s}+\Delta \lambda \hat \Gamma'  \right )        
\end{equation}

\begin{equation}
\der{f}{\theta_s}=\left(\Gamma' q +  \Gamma \der{q }{\theta_s}\right)+M  \der{p}{\theta_s}-\der{\sigma_0}{\Delta \lambda}\der{\Delta \lambda}{\theta_s}
\end{equation}

\section{Consistent Jacobian operators}

\subsubsection{Consistent operator for the linear elastic response}
In case of linear elastic response, the stresses $\tenss{\sigma}$  and couple-stresses $\tenss{\mu}$ at the end of the increment 
 coincide with their respective  elastic predictors $\tenss{\sigma}^*$ and $\tenss{\mu}^*$ and the consistent tangent operators coincide with  the elastic fourth order stiffness tensors.
 
%

\subsubsection{Consistent operator for the first return algorithm}
 By differentiating Eqs. \eqref{Eq_StressTensorsCosseratSymFirstAlgo}
\begin{equation}
\begin{gathered}
\begin{aligned}
\der{\tenss{\sigma}_{sym}}{\tenss{\varepsilon} }=\frac{1}{q^*} \left(\tenss{s}_{sym}^* \otimes \der{q}{\tenss{\varepsilon}} \right)-
3 G \frac{q}{{q^*}^3} \left( \tenss{s}_{sym}^*  \otimes \tenss{s}_{sym}^*  \right) + 2G\frac{q }{q^*}\tensff{ \mathcal{I}}_d
+\left( \tenss{I}  \otimes \der{p}{{\tenss{\varepsilon} }}\right)
\end{aligned}\\
\der{\tenss{\sigma}_{sym}}{\tenss{\omega} }=\frac{1}{q^*} \left(\tenss{s}_{sym}^* \otimes \der{q}{\tenss{\omega}}\right)-
3G  \frac{q}{{q^*}^3} \left(\tenss{s}_{sym}^* \otimes \tenss{s}_{skw}^*\right) + \left(\tenss{I} \otimes \der{p}{\tenss{\omega}}\right)\\
\der{\tenss{\sigma}_{sym}}{\tenss{\chi}}=\frac{1}{q^*} \left(\tenss{s}_{sym}^* \otimes \der{q}{\tenss{\chi}}\right)-
3G  \frac{q}{{q^*}^3} \left(\tenss{s}_{sym}^*  \otimes \tenss{\mu}^*\right) + \left( \tenss{I}  \otimes \der{p}{\tenss{\chi}}\right)\\
\end{gathered}
\end{equation}
and  by differentiating  Eq. \eqref{Eq_StressTensorsCosseratComponents}, using the second, third and forth of Eq. \eqref{Eq_StressTensorsFirstAlghorithm}
\begin{equation}
\begin{gathered}
\der{\tenss{s}_{skw}}{\tenss{\varepsilon}}=\frac{1}{q^*} \left(\tenss{s}_{skw}^* \otimes \der{q}{{\tenss{\varepsilon}}}\right)-
3G  \frac{q}{{q^*}^3} \left(\tenss{s}_{skw}^* \otimes \tenss{s}_{sym}^*\right)\\
\der{ \tenss{s}_{skw} }{\tenss{\omega} }=\frac{1}{q^*} \left( \tenss{s}_{skw}^* \otimes \der{q}{\tenss{\omega} }\right)-
3G \frac{q}{{q^*}^3} \left(\tenss{s}_{skw}^* \otimes \tenss{s}_{skw}^*\right)+2 G_c\frac{q}{q^*} \tensff{I} \\
\der{\tenss{s}_{skw}}{\tenss{\chi}}=\frac{1}{q^*} \left(\tenss{s}_{skw}^* \otimes \der{q}{\tenss{\chi}}\right)-
3G \frac{q}{{q^*}^3} \left(\tenss{s}_{skw}^* \otimes \tenss{\mu}^* \right)\\
\end{gathered}
\end{equation}

\begin{equation}
\begin{gathered}
\der{\tenss{\mu}}{\tenss{\varepsilon} }=\frac{1}{q^*} \left(\tenss{\mu}^* \otimes \der{q}{\tenss{\varepsilon}}\right)-
3G \frac{q}{{q^*}^3} \left(\tenss{\mu}^* \otimes \tenss{s}_{sym}^*\right)\\
\der{\tenss{\mu}}{\tenss{\omega}}=\frac{1}{q^*} \left(\tenss{\mu}^* \otimes \der{q}{\tenss{\omega}}\right)-
3G \frac{q}{{q^*}^3} \left(\tenss{\mu}^* \otimes \tenss{s}_{skw}^*\right)\\
\begin{aligned}
\der{\tenss{\mu}}{\tenss{\chi}}=\frac{1}{q^*} \left(\tenss{\mu}^* \otimes \der{q}{\tenss{\chi}}\right)-
3G  \frac{q}{{q^*}^3} \left(\tenss{\mu}^* \otimes \tenss{\mu}^*\right)+
\frac{q}{q^*} \left[2B \left( \bar{\tensff{I}}^{sym} - \frac{1}{3} \bar{\bar{\tensff{I}}}  \right)  + 2 B_c \bar{\tensff{I} }^{skw} \right]+ K_c  \bar{\bar{\tensff{I}}} 
\end{aligned}
\\
\end{gathered}
\end{equation}

Where the required derivatives are 

\begin{equation}
\begin{gathered}
\der{p(\tenss{\varepsilon}, \tenss{\omega}, \tenss{\chi})}{\tenss{\varepsilon} }=
\der{p\left(\tenss{\varepsilon}, \Delta \lambda (\tenss{\varepsilon}, \tenss{\omega}, \tenss{\chi}) \right)}{\tenss{\varepsilon} }=
\der{p}{\tenss{\varepsilon}}+\der{p}{\Delta \lambda} \der{\Delta \lambda}{{\tenss{\varepsilon}}}=K \left( \tenss{ I} - \hat M \der{\Delta \lambda}{\tenss{\varepsilon}}\right)\\
\der{p(\tenss{\varepsilon}, \tenss{\omega}, \tenss{\chi})}{\tenss{\omega} }=
\der{p\left(\tenss{\varepsilon}, \Delta \lambda (\tenss{\varepsilon}, \tenss{\omega}, \tenss{\chi}) \right)}{\tenss{\omega} }=
\der{p}{\Delta \lambda} \der{\Delta \lambda}{{\tenss{\omega}}}= -K  \hat M \der{\Delta \lambda}{\tenss{\omega}}\\
\der{p(\tenss{\varepsilon}, \tenss{\omega}, \tenss{\chi})}{\tenss{\chi} }=
\der{p\left(\tenss{\varepsilon}, \Delta \lambda (\tenss{\varepsilon}, \tenss{\omega}, \tenss{\chi}) \right)}{\tenss{\chi} }=
\der{p}{\Delta \lambda} \der{\Delta \lambda}{{\tenss{\chi}}}=- K  \hat M \der{\Delta \lambda}{\tenss{\chi}}\\
\end{gathered}
\end{equation}

\begin{equation}
\begin{gathered}
\der{q \left(\tenss{\varepsilon}, \tenss{\omega}, \tenss{\chi} \right) }{\tenss{\varepsilon} }=
\der{q \left(\tenss{\varepsilon}, \tenss{\omega}, \tenss{\chi}, \Delta\lambda(\tenss{\varepsilon}, \tenss{\omega}, \tenss{\chi}) \right) }{\tenss{\varepsilon} }=
\der{q}{\tenss{\varepsilon}}+\der{q}{\Delta \lambda} \der{\Delta \lambda}{{\tenss{\varepsilon}}}=\\
 3G \left( \frac{\tenss{s}_{sym}^*}{q^*}- \hat \Gamma\left(\theta^*_s\right)\der{\Delta \lambda}{{\tenss{\varepsilon}}}\right)\\
\der{q\left(\tenss{\varepsilon}, \tenss{\omega}, \tenss{\chi} \right) }{\tenss{\omega}}=
\der{q \left(\tenss{\varepsilon}, \tenss{\omega}, \tenss{\chi}, \Delta\lambda(\tenss{\varepsilon}, \tenss{\omega}, \tenss{\chi}) \right) }{\tenss{\omega}}=
\der{q}{\tenss{\omega}}+\der{q}{\Delta \lambda} \der{\Delta \lambda}{\tenss{\omega}}=\\
3G  \left( \frac{\tenss{s}_{skw}^*}{q^*}- \hat \Gamma\left(\theta^*_s\right)\der{\Delta \lambda}{\tenss{\omega}}\right)\\
\der{q  \left(\tenss{\varepsilon}, \tenss{\omega}, \tenss{\chi} \right)}{\tenss{\chi}}=
\der{q \left(\tenss{\varepsilon}, \tenss{\omega}, \tenss{\chi}, \Delta\lambda(\tenss{\varepsilon}, \tenss{\omega}, \tenss{\chi}) \right)  }{\tenss{\chi}}=
\der{q}{\tenss{\chi}}+\der{q}{\Delta \lambda} \der{\Delta \lambda}{\tenss{\chi}}= \\
3G  \left( \frac{\tenss{\mu}^*}{q^*}- \hat \Gamma\left(\theta^*_s\right)\der{\Delta \lambda}{\tenss{\chi}}\right)\\
\end{gathered}
\end{equation}

\begin{equation}
\der{p}{\Delta \lambda} = -K \hat M
\end{equation}

\begin{equation}
\der{q}{\Delta \lambda} =-3G  \hat \Gamma\left(\theta^*_s\right)
\end{equation}

\begin{equation}
\der{\Delta \lambda}{{\tenss{\varepsilon}}}= -\frac{\der{f}{\tenss{\varepsilon}}}{\der{f}{\Delta \lambda}}
=
\frac{\ds 3G  \ds \frac{\tenss{s}_{sym}^*}{q^*}  \Gamma \left(\theta^*_s \right)+ M K \tenss{ I} }{\ds 3G \hat \Gamma\left(\theta^*_s\right)    \Gamma \left(\theta^*_s \right)  + K M \hat M  +\der{\sigma_0}{\Delta \lambda} }
\end{equation}

\begin{equation}
\der{\Delta \lambda}{{\tenss{\omega}}}= -\frac{\der{f}{\tenss{\omega}}}{\der{f}{\Delta \lambda}}
=
\frac{\ds 3G  \ds \frac{\tenss{s}_{skw}^*}{q^*} \Gamma \left(\theta^*_s \right)}{\ds 3G \hat \Gamma\left(\theta^*_s\right)    \Gamma \left(\theta^*_s \right)  + K M \hat M  +\der{\sigma_0}{\Delta \lambda} }
\end{equation}
\begin{equation}
\der{\Delta \lambda}{\tenss{\chi}}= -\frac{\der{f}{\tenss{\chi}}}{\der{f}{\Delta \lambda}}
=\frac{\ds 3G  \ds \frac{\tenss{\mu}^*}{q^*}  \Gamma \left(\theta^*_s \right)}{\ds 3G \hat \Gamma\left(\theta^*_s\right)    \Gamma \left(\theta^*_s \right)  + K M \hat M  +\der{\sigma_0}{\Delta \lambda} }
\end{equation}

\begin{equation}
\der{f}{\Delta \lambda}=  -3G \hat \Gamma\left(\theta^*_s\right)    \Gamma \left(\theta^*_s \right)   - K M \hat M  -\der{\sigma_0}{\Delta \lambda} 
\end{equation}

\subsubsection{Consistent operator for the second return algorithm }
The construction of the consistent tangent operator in the general case of the second return algorithm requires the evaluation of the derivative of  Eq. \eqref{Eq_spectral} with respect to $\tenss{\varepsilon}$ 
\begin{equation}
\begin{aligned}
\der{\tenss{\sigma}_{sym}}{\tenss{\varepsilon}}=\sum_{i=\rn{1}, \rn{2}, \rn{3}} \left(\tenss{b}_i^* \otimes \der{\sigma_i}{\tenss{\varepsilon}}\right) +\sigma_i \tensff{\Omega}_i^*
\end{aligned}
\end{equation}
where $\tensff{\Omega}_i^*= \derl{\tenss{b}_i^*}{\tenss{\varepsilon}}$ is the spin of the principal directions. 
Since, at variance with standard approaches, all principal stresses in the general case of the second return algorithm are distinct,  and the principal directions of the symmetric stress state $\tenss{\sigma}_{sym}$  at time $t_{n+1}$ coincide with those of its predictor $\tenss{\sigma}_{sym}^*$ or, equivalently, of $\tenss{\varepsilon}^*$,  the evaluation of $\tensff{\Omega}_i^*$ following Miehe  \cite{Miehe1998} is simply 
\begin{equation}
\tensff{\Omega}^*_i=\der{\tenss {b}^*_i}{\tenss{\varepsilon}}= \der{\tenss {b}^*_i}{\tenss{\varepsilon}^*}\equiv \der{\varepsilon^*_i}{\tenss{\varepsilon}^*_s \otimes \tenss{\varepsilon}^*_s}
\end{equation}
The derivatives of $\tenss{\sigma}_{sym}$ with respect the other Cosserat strain tensors is similarly obtained from  Eq. \eqref{Eq_spectral}
\begin{equation}
\begin{gathered}
\der{\tenss{\sigma}_{sym}}{\tenss{\omega}}=\sum_{i=\rn{1}, \rn{2}, \rn{3}} \left(\tenss{b}^*_i \otimes \der{\sigma_i}{\tenss{\omega}}\right) \\
\der{\tenss{\sigma}_{sym}}{\tenss{\chi}}=\sum_{i=\rn{1}, \rn{2}, \rn{3}} \left(\tenss{b}^*_i \otimes \der{\sigma_i}{\tenss{\chi}}\right)
\end{gathered}
\end{equation}
The others consistent jacobians operators are constructed by deriving  the third, fourth and fifth of Eqs. \eqref{Eq_StressesFuncPredicLambdaTheta} with respect to each of the Cosserat strain tensors
\begin{equation}
\begin{gathered}
\der{\tenss{s}_{skw}}{\tenss{\varepsilon} }=\tenss{s}_{skw}^* \otimes  \left(\frac{1}{r}  \der{q}{\tenss{\varepsilon}}  -\frac{q}{r^2}  \der{r}{\tenss{\varepsilon} }\right)\\
\der{\tenss{s}_{skw} }{\tenss{\omega}}=\tenss{s}_{skw}^* \otimes  \left(\frac{1}{r}  \der{q}{\tenss{\omega}}  -\frac{q}{r^2}  \der{r}{\tenss{\omega}}\right) +2 G_c \frac{q}{r} \tensff{I}  \\
\der{\tenss{s}_{skw}}{\tenss{\chi}}=\tenss{s}_{skw}^* \otimes  \left(\frac{1}{r}  \der{q}{\tenss{\chi}}  -\frac{q}{r^2}  \der{r}{\tenss{\chi}}\right)\\
\der{\tenss{\mu}}{\tenss{\varepsilon} }=\tenss{\mu}^* \otimes  \left(\frac{1}{r}  \der{q}{\tenss{\varepsilon}}  -\frac{q}{r^2}  \der{r}{\tenss{\varepsilon}}\right) \\
\der{\tenss{\mu}}{\tenss{\omega}}=\tenss{\mu}^* \otimes  \left(\frac{1}{r}  \der{q}{\tenss{\omega}}  -\frac{q}{r^2}  \der{r}{\tenss{\omega}}\right)\\
\begin{aligned}
\der{\tenss{\mu}}{\tenss{\chi}}=\tenss{\mu}^* \otimes  \left(\frac{1}{r}  \der{q}{\tenss{\chi}}  -\frac{q}{r^2}  \der{r}{\tenss{\chi}}\right)+ \frac{q}{r} \left[2B \left( \bar{\tensff{I}}^{sym} - \frac{1}{3} \bar{\bar{\tensff{I}}}  \right)  + 2 B_c \bar{\tensff{I} }^{skw} \right]+ K_c  \bar{\bar{\tensff{I}}} 
\end{aligned}
\end{gathered}
\end{equation}
whilst all required derivatives are 

Let consider that the three invaria
\begin{equation}
\begin{gathered}
\der{r(\theta_s,\tenss{\varepsilon}, \tenss{\omega}, \tenss{\chi})}{\tenss{\varepsilon}}=3G \cos^2 \left(\theta_s^*-\theta_s\right)\frac{\tenss{s}_{sym}^*}{r}-\der{r}{\theta_s}\der{\theta_s^*}{\tenss{\varepsilon}}\\
\der{r(\theta_s,\tenss{\varepsilon}, \tenss{\omega}, \tenss{\chi})}{\tenss{\omega}}=3G \frac{\tenss{s}_{skw}^*}{r}\\
\der{r(\theta_s,\tenss{\varepsilon}, \tenss{\omega}, \tenss{\chi})}{\tenss{\chi}}=3G \frac{\tenss{\mu}^*}{r}\\
\end{gathered}
\end{equation}

\begin{equation}
\begin{gathered}
\der{\Delta \lambda}{\tenss{\varepsilon}}= \Delta \lambda \left[ 6 G \frac{\tenss{s}_{sym}^*}{{q_s^*}^2}- \frac{1}{r}\der{r(\theta_s,\tenss{\varepsilon}, \tenss{\omega}, \tenss{\chi})}{\tenss{\varepsilon}} +
\frac{2}{\tan\left(2\theta_s^*-2\theta_s\right)}\der{\theta_s^*}{\tenss{\varepsilon} }\right]\\
\der{\Delta \lambda}{\tenss{\omega}}=-\frac{\Delta \lambda}{r} \der{r(\theta_s,\tenss{\varepsilon}, \tenss{\omega}, \tenss{\chi})}{\tenss{\omega} }\\
\der{\Delta \lambda}{\tenss{\chi}}=-\frac{\Delta \lambda}{r} \der{r(\theta_s,\tenss{\varepsilon}, \tenss{\omega}, \tenss{\chi})}{\tenss{\chi}}
\end{gathered}
\end{equation}

\begin{equation}
\begin{gathered}
\der{p(\theta_s,\tenss{\varepsilon}, \tenss{\omega}, \tenss{\chi})}{\tenss{\varepsilon}}=K \left(\tensf{I} 
-\hat M \der{\Delta \lambda}{\tenss{\varepsilon}}\right)\\
\der{p(\theta_s,\tenss{\varepsilon}, \tenss{\omega}, \tenss{\chi})}{\tenss{\omega}}=- K \hat M \der{\Delta \lambda}{\tenss{\omega}}\\
\der{p(\theta_s,\tenss{\varepsilon}, \tenss{\omega}, \tenss{\chi})}{\tenss{\chi}}=- K \hat M \der{\Delta \lambda}{\tenss{\chi}}
\end{gathered}
\end{equation}
\begin{equation}
\begin{gathered}
\der{q(\theta_s,\tenss{\varepsilon}, \tenss{\omega}, \tenss{\chi})}{\tenss{\varepsilon}}=\der{r(\theta_s,\tenss{\varepsilon}, \tenss{\omega}, \tenss{\chi})}{\tenss{\varepsilon}}-3G  \hat \Gamma  \der{\Delta \lambda}{\tenss{\varepsilon}}\\
\der{q(\theta_s,\tenss{\varepsilon}, \tenss{\omega}, \tenss{\chi})}{\tenss{\omega}}=\der{r(\theta_s,\tenss{\varepsilon}, \tenss{\omega}, \tenss{\chi})}{\tenss{\omega}}-3G  \hat \Gamma  \der{\Delta \lambda}{\tenss{\omega}}\\
\der{q(\theta_s,\tenss{\varepsilon}, \tenss{\omega}, \tenss{\chi})}{\tenss{\chi}}=\der{r(\theta_s,\tenss{\varepsilon}, \tenss{\omega}, \tenss{\chi})}{\tenss{\chi}}-3G   \hat \Gamma  \der{\Delta \lambda}{\tenss{\chi}}\\
\end{gathered}
\end{equation}

\begin{equation}
\begin{gathered}
\begin{aligned}
\der{\theta_s}{\tenss{\varepsilon}}=-\frac{\der{f}{\tepsilon_s} }{\der{f}{\theta_s} }
=-
\frac{1}{\der{f}{\theta_s} }\left( \Gamma \der{q(\theta_s,\tenss{\varepsilon}, \tenss{\omega}, \tenss{\chi})}{\tenss{\varepsilon}}+M \der{p(\theta_s,\tenss{\varepsilon}, \tenss{\omega}, \tenss{\chi})}{\tenss{\varepsilon}}
\right. \\ \left.
-\der{\sigma_0}{\Delta \lambda}
\der{\Delta \lambda}{\tenss{\varepsilon}}\right)
\end{aligned}
\\
\begin{aligned}
\der{\theta_s}{\tenss{\omega}}=-\frac{\der{f}{\tenss{\omega}} }{\der{f}{\theta_s} }
=-
\frac{1}{\der{f}{\theta_s} }\left( \Gamma \der{q(\theta_s,\tenss{\varepsilon}, \tenss{\omega}, \tenss{\chi})}{\tenss{\omega}}+M \der{p(\theta_s,\tenss{\varepsilon}, \tenss{\omega}, \tenss{\chi})}{\tenss{\omega}}
\right. \\ \left.-\der{\sigma_0}{\Delta \lambda}
\der{\Delta \lambda}{\tenss{\omega}}\right)
\end{aligned}
\\
\begin{aligned}
\der{\theta_s}{\tenss{\chi}}=-\frac{\der{f}{\tenss{\chi}} }{\der{f}{\theta_s} }
=-
\frac{1}{\der{f}{\theta_s} }\left( \Gamma \der{q(\theta_s,\tenss{\varepsilon}, \tenss{\omega}, \tenss{\chi})}{\tenss{\chi}}+M \der{p(\theta_s,\tenss{\varepsilon}, \tenss{\omega}, \tenss{\chi})}{\tenss{\chi}}
\right. \\ \left.-\der{\sigma_0}{\Delta \lambda}
\der{\Delta \lambda}{\tenss{\chi}}\right)
\end{aligned}
\end{gathered}
\end{equation}

\begin{equation}
\begin{gathered}
\der{p(\tenss{\varepsilon}, \tenss{\omega}, \tenss{\chi})}{\tenss{\varepsilon}}=
\der{p\left(\tenss{\varepsilon}, \tenss{\omega}, \tenss{\chi}, \theta_s(\tenss{\varepsilon}, \tenss{\omega}, \tenss{\chi}) \right)}{\tenss{\varepsilon}}=
\der{p}{\tenss{\varepsilon}}+\der{p}{\theta_s}\der{\theta_s}{\tenss{\varepsilon}}\\
\der{p(\tenss{\varepsilon}, \tenss{\omega}, \tenss{\chi})}{\tenss{\omega}}=
\der{p\left(\tenss{\varepsilon}, \tenss{\omega}, \tenss{\chi}, \theta_s(\tenss{\varepsilon}, \tenss{\omega}, \tenss{\chi}) \right)}{\tenss{\omega}}=
\der{p}{\tenss{\omega}}+\der{p}{\theta_s}\der{\theta_s}{\tenss{\omega}}\\
\der{p(\tenss{\varepsilon}, \tenss{\omega}, \tenss{\chi})}{\tenss{\chi}}=
\der{p\left(\tenss{\varepsilon}, \tenss{\omega}, \tenss{\chi}, \theta_s(\tenss{\varepsilon}, \tenss{\omega}, \tenss{\chi}) \right)}{\tenss{\chi}}=
\der{p}{\tenss{\chi}}+\der{p}{\theta_s}\der{\theta_s}{\tenss{\chi}}\\
\end{gathered}
\end{equation}

\begin{equation}
\begin{gathered}
\der{q(\tenss{\varepsilon}, \tenss{\omega}, \tenss{\chi})}{\tenss{\varepsilon}}=
\der{q\left(\tenss{\varepsilon}, \tenss{\omega}, \tenss{\chi}, \theta_s(\tenss{\varepsilon}, \tenss{\omega}, \tenss{\chi}) \right)}{\tenss{\varepsilon}}=
\der{q}{\tenss{\varepsilon}}+\der{q}{\theta_s}\der{\theta_s}{\tenss{\varepsilon}}\\
\der{q(\tenss{\varepsilon}, \tenss{\omega}, \tenss{\chi})}{\tenss{\omega}}=
\der{q\left(\tenss{\varepsilon}, \tenss{\omega}, \tenss{\chi}, \theta_s(\tenss{\varepsilon}, \tenss{\omega}, \tenss{\chi}) \right)}{\tenss{\omega}}=
\der{q}{\tenss{\omega}}+\der{q}{\theta_s}\der{\theta_s}{\tenss{\omega}}\\
\der{q(\tenss{\varepsilon}, \tenss{\omega}, \tenss{\chi})}{\tenss{\chi}}=
\der{q\left(\tenss{\varepsilon}, \tenss{\omega}, \tenss{\chi}, \theta_s(\tenss{\varepsilon}, \tenss{\omega}, \tenss{\chi}) \right)}{\tenss{\chi}}=
\der{q}{\tenss{\chi}}+\der{q}{\theta_s}\der{\theta_s}{\tenss{\chi}}\\
\end{gathered}
\end{equation}

\begin{equation}
\begin{gathered}
\der{r(\tenss{\varepsilon}, \tenss{\omega}, \tenss{\chi})}{\tenss{\varepsilon}}=
\der{r\left(\tenss{\varepsilon}, \tenss{\omega}, \tenss{\chi}, \theta_s(\tenss{\varepsilon}, \tenss{\omega}, \tenss{\chi}) \right)}{\tenss{\varepsilon}}=
\der{r}{\tenss{\varepsilon}}+\der{r}{\theta_s}\der{\theta_s}{\tenss{\varepsilon}}\\
\der{r(\tenss{\varepsilon}, \tenss{\omega}, \tenss{\chi})}{\tenss{\omega}}=
\der{r\left(\tenss{\varepsilon}, \tenss{\omega}, \tenss{\chi}, \theta_s(\tenss{\varepsilon}, \tenss{\omega}, \tenss{\chi}) \right)}{\tenss{\omega}}=
\der{r}{\tenss{\omega}}+\der{r}{\theta_s}\der{\theta_s}{\tenss{\omega}}\\
\der{r(\tenss{\varepsilon}, \tenss{\omega}, \tenss{\chi})}{\tenss{\chi}}=
\der{r\left(\tenss{\varepsilon}, \tenss{\omega}, \tenss{\chi}, \theta_s(\tenss{\varepsilon}, \tenss{\omega}, \tenss{\chi}) \right)}{\tenss{\chi}}=
\der{r}{\tenss{\chi}}+\der{r}{\theta_s}\der{\theta_s}{\tenss{\chi}}\\
\end{gathered}
\end{equation}

\begin{equation}
\begin{gathered}
\begin{aligned}
\der{q_s}{\tenss{\varepsilon}}=\frac{\cos\left(\theta_s^*-\theta_s\right)}{r}
\left[
q_s^*\der{q(\tenss{\varepsilon}, \tenss{\omega}, \tenss{\chi})}{\tenss{\varepsilon}}+ 3G q
\frac{\tenss{s}_{sym}^*}{q_s^*}
\right]\\
+
\frac{q q_s^* \sin\left(\theta_s^*-\theta_s\right)}{r}\left(\der{\theta_s}{\tenss{\varepsilon}}-\der{\theta_s^*}{\tenss{\varepsilon}} \right)-\frac{q_s}{r}\der{r(\tenss{\varepsilon}, \tenss{\omega}, \tenss{\chi})}{\tenss{\varepsilon}}
\end{aligned}\\
\begin{aligned}
\der{q_s}{\tenss{\omega}}=\frac{q_s^* \cos\left(\theta_s^*-\theta_s\right)}{r}\der{q(\tenss{\varepsilon}, \tenss{\omega}, \tenss{\chi})}{\tenss{\omega}}+
\frac{q q_s^* \sin\left(\theta_s^*-\theta_s\right)}{r}\der{\theta_s}{\tenss{\omega}}\\-\frac{q_s}{r}\der{r(\tenss{\varepsilon}, \tenss{\omega},\tenss{\chi} )}{\tenss{\omega}}
\end{aligned}\\
\begin{aligned}
\der{q_s}{\tenss{\chi}}=\frac{q_s^* \cos\left(\theta_s^*-\theta_s\right)}{r}\der{q(\tenss{\varepsilon}, \tenss{\omega} , \tenss{\chi})}{\tenss{\chi}}+
\frac{q q_s^* \sin\left(\theta_s^*-\theta_s\right)}{r}\der{\theta_s}{\tenss{\chi}}\\-\frac{q_s}{r}\der{r(\tenss{\varepsilon}, \tenss{\omega},\tenss{\chi})}{\tenss{\chi}}
\end{aligned}\\
\end{gathered}
\end{equation}

\begin{equation}
\begin{aligned}
\begin{gathered}
\der{\sigma_i}{\tenss{\varepsilon}}=\der{p(\tenss{\varepsilon}, \tenss{\omega}, \tenss{\chi} )}{\tenss{\varepsilon}}+
\frac{2}{3}\left[\sin\left(\beta_i \right)\der{q_s}{\tenss{\varepsilon}}+q_s
\cos \left(\beta_i \right) \der{\theta_s} {\tenss{\varepsilon}}\right]\\
\der{\sigma_i}{\tenss{\omega}}=\der{p(\tenss{\varepsilon}, \tenss{\omega}, \tenss{\chi} )}{\tenss{\omega}}+
\frac{2}{3}\left[\sin\left(\beta_i \right)\der{q_s}{\tenss{\omega}}+q_s
\cos \left(\beta_i \right) \der{\theta_s} {\tenss{\omega}}\right]\\
\der{\sigma_i}{\tenss{\chi}}=\der{p(\tenss{\varepsilon}, \tenss{\omega}, \tenss{\chi} )}{\tenss{\chi}}+
\frac{2}{3}\left[\sin\left(\beta_i \right)\der{q_s}{\tenss{\chi}}+q_s
\cos \left(\beta_i \right) \der{\theta_s} {\tenss{\chi}}\right]\\
\end{gathered}\\
\beta_\rn{1}=\left(\theta_s+\frac{2}{3}\pi \right), \;\;
\beta_\rn{2}=\theta_s, \;\;
\beta_\rn{3}=\left(\theta_s-\frac{2}{3}\pi \right)
\end{aligned}
\end{equation}

The computation of the derivatives of the stress invariants $p$, $q$, $\theta_s$ and of the principal stress components with respect to the strain increments is standard, and it is reported in Appendix.

\subsubsection{Consistent operator for the third return algorithm}
\begin{equation}
\der{\tenss{\sigma}_{sym}}{\tenss{\varepsilon}}= K \left[ 1- \left(\frac{ M \hat M K }{K M \hat M+ \der{\sigma_0}{\Delta \lambda}} \right)\right] \bar{\bar{ \tensff{I} }}
\end{equation}

\begin{equation}
\der{\tenss{\sigma}_{sym}}{\tenss{\omega}}=\der{\tenss{\sigma}_{sym}}{\tenss{\chi}}=\tensff{0},\;\;\;
\der{\tenss{s}_{skw}}{\tenss{\varepsilon}}=\der{\tenss{s}_{skw}}
{\tenss{\omega}}=\der{\tenss{s}_{skw}}{\tenss{\chi}}=\tensff{0},\;\;\;
\der{\tenss{\mu}}{\tenss{\varepsilon}}=\der{\tenss{\mu}}{\tenss{\omega}}=\tensff{0}
\end{equation}
\begin{equation}
\der{\tenss{\mu}}{\tenss{\chi}}= \tensff{0} 
\end{equation}

\section {Matrices for 2D plane strain FE}
By ordering the generalized nodal displacements components in the vector
\begin{equation}
\tensf{u}_g= [\underbrace{u_x^{(1)}\;\; u_y^{(1)}\;\; \theta_z^{(1)}\;\;...\;\; u_x^{(\hat n)}\;\; u_y^{(\hat n)}\; \theta_z^{(\hat n)}}_{3 \hat n} ]^T
\end{equation}
one can compute the vectors $\tensf{\varepsilon} $ and $\tensf{\omega}$ containing the components of the tensors $\tenss{\varepsilon}$ and $\tenss{\omega}$ respectively as:
\begin{equation}
\tensf{\varepsilon}= \left[ \der{u_x}{x}\;\; \der{u_y}{y}\;\;0 \;\; \frac{1}{2} \left( \der{u_y}{x}+ \der {u_x}{y} \right)\;\;\frac{1}{2} \left(  \der {u_x}{y} + \der{u_y}{x} \right) \right]^T= \bvec B \hat{\tensf{ u}}_g
\end{equation}
\begin{equation}
\tensf{\omega} = \left[0\;\;0\;\;0\;\; \frac{1}{2}\left( \der{u_y}{x} - \der{u_x}{y} \right)-\theta_z\;\;    \frac{1}{2}\left( \der{u_x}{y} - \der{u_y}{x} \right)+\theta_z \right]^T= \bvec W \hat{\tensf{ u}}_g
\end{equation}
where:
\begin{equation}
\bvec B= 
\begin{bmatrix}
\der{N^{(1)}}{x}             & 0                           & 0 & .. & \der{N^{(\hat n)}}{x}            &  0                                & 0 \\
0                            & \der{N^{(1)}}{y}            & 0 & .. & 0                                & \der{N^{(\hat n)}}{y}             & 0 \\
0                            & 0                           & 0 & .. & 0                                &  0                                & 0 \\
\ds\frac{1}{2} \der{N^{(1)}}{y} & \ds\frac{1}{2}\der{N^{(1)}}{x} & 0 & .. &\ds\frac{1}{2} \der{N^{(\hat n)}}{y} & \ds\frac{1}{2} \der{N^{(\hat n)}}{x} & 0\\
\ds\ds\frac{1}{2} \der{N^{(1)}}{y} & \ds\frac{1}{2}\der{N^{(1)}}{x} & 0 & .. &\ds\frac{1}{2} \der{N^{(\hat n)}}{y} & \ds\frac{1}{2} \der{N^{(\hat n)}}{x} & 0
\end{bmatrix}
\end{equation}
\begin{equation}
\bvec W= 
\begin{bmatrix}
0                            & 0                           & 0 & .. & 0                                &  0                                & 0 \\
0                            & 0                           & 0 & .. & 0                                &  0                                & 0 \\
0                            & 0                           & 0 & .. & 0                                &  0                                & 0 \\
-\ds\frac{1}{2} \der{N^{(1)}}{y} & \ds\frac{1}{2}\der{N^{(1)}}{x} & -N^{(1)} & .. &-\ds\frac{1}{2} \der{N^{(\hat n)}}{y} & \ds\frac{1}{2} \der{N^{(\hat n)}}{x} & -N^{(\hat n)}\\
\ds\ds\frac{1}{2} \der{N^{(1)}}{y} & -\ds\frac{1}{2}\der{N^{(1)}}{x} & N^{(1)} & .. &\ds\frac{1}{2} \der{N^{(\hat n)}}{y} & -\ds\frac{1}{2} \der{N^{(\hat n)}}{x} & N^{(\hat n)}
\end{bmatrix}
\end{equation}
The vector $\tensf{\chi}$ containing the components of the tensor $\tenss{\chi}$ can be computed as:
\begin{equation}
\tensf{\chi}= \left[\der{\theta_z}{x}\;\;\der{\theta_z}{y}\right]^T= \bvec M \tensf{u}_g
\end{equation}
\begin{equation}
\bvec M= 
\begin{bmatrix}
0 &0 &  \ds \der{N^{(1)}}{x}  & .. & 0 & 0 & \ds \der{N^{(\hat n)}}{x}\\
0 &0 &  \ds \der{N^{(1)}}{y}  & .. & 0 & 0 & \ds \der{N^{(\hat n)}}{y}\\
\end{bmatrix}
\end{equation}


%

\end{appendices}

\end{document}